\documentclass[
11pt, oneside, 
a4paper]{article}
\usepackage[ansinew]{inputenc}
\usepackage[english]{babel}
\usepackage{comment,amsbsy,amscd,amsfonts,amsmath,amsopn,amssymb,amstext,amsthm,amsxtra,array,color}
\usepackage{fancybox,float,graphicx,latexsym,subfigure,srcltx,times,tikz,url,amsthm}
\usepackage{bm}
\newcommand{\mx}[1]{W_n}

\newcommand{\ver}{\operatorname{vert}}

\newcommand{\N}{\mathbb{N}}
\newcommand{\Z}{\mathbb{Z}}
\newcommand{\Q}{\mathbb{Q}}
\newcommand{\R}{\mathbb{R}}

\newcommand{\FF}{\mathcal{F}}

\newcommand{\II}{\mathcal{I}}
\newcommand{\JJ}{\mathcal{J}}

\newcommand{\MM}{\mathcal{M}}

\newcommand{\PP}{\mathcal{P}}

\newcommand{\se}{\operatorname{secc}}

\newcommand{\Roof}{\operatorname{Roof}}

\newcommand{\rad}{\operatorname{rad}}
\newcommand{\pyr}{\operatorname{pyr}}

\newcommand{\aff}{\operatorname{aff}}

\newcommand{\dd}{\operatorname{d}}

\newcommand{\secc}{\operatorname{secc}}

\newcommand{\conv}{\operatorname{conv}}

\newcommand{\diag}{\operatorname{diag}}

\newcommand{\vol}{\operatorname{vol}}

\newtheorem{thm}{Theorem}[section]
\newtheorem{lem}[thm]{Lemma}
\newtheorem{cor}[thm]{Corollary}
\newtheorem{prop}[thm]{Proposition}
\newtheorem{ex}[thm]{Example}

\newtheorem{rem}[thm]{Remark}

\theoremstyle{definition}
\newtheorem{dfn}[thm]{Definition}
\newtheorem{nota}[thm]{Notation}

\newtheorem{question}{Question}

\graphicspath{{images/}{C:/Users/user/Documents/docus/tropical_nueva/todas_las_figuras_pdf/las_que_uso_en_volume/}
{C:/Users/user/Documents/docus/tropical_nueva/todas_las_figuras_jpg/las_que_uso_en_volume/}}

\title{The volume of an isocanted   cube is a determinant\footnote{In the talk given in ALAMA 2022, \emph{isocanted cubes} were given the longer name of  \emph{isocanted alcoved polytopes.} }}

\author{M.J. de la Puente\thanks{Partially supported by Ministerio de Ciencia e Innovaci\'{o}n, PID2019-107701GB-I00 and by UCM,
Research group 910444, and by RED 2022-134176-T.} \thanks{Corresponding author.}\\ Dpto. de \'{A}lgebra, Geometr\'{\i}a y Topolog\'{\i}a\\
Facultad de Matem\'{a}ticas\\Universidad Complutense (Spain)\\\texttt{mpuente@ucm.es}\\ Phone: 34--91--3944659\\
P.L. Claver\'{i}a\\Universidad de Zaragoza (Spain)\\\texttt{plcv@unizar.es}}

\begin{document}
\maketitle
\begin{abstract}
In any dimension $d\ge2$, we give exact  volume formulas of two mutually polar dual convex $d$--polytopes. The primal body
is  called \emph{isocanted cube} of dimension $d$, depending on
two \emph{real} parameters $0<a<\ell$. The limit case $a=0$ yields a $d$--cube of edge--length  $\ell$. We prove that the  volume of
such a body is
 the determinant of the matrix  of order $d$ having   diagonal  entries equal to
$\ell$ and  $a$ elsewhere.

We also compute
the volume of the   polar dual body,  getting a rational expression in $\ell,a$, homogeneous of degree $-d$.
with rational coefficients.

Isocanted cubes are origin--symmetric
zonotopes.
Zonoids (defined as  the limits of families of
zonotopes) satisfy  \emph{the Mahler conjecture}; in particular,  zonotopes do.
Nonetheless, we confirm (by elementary methods) that  \emph{the Mahler conjecture}
holds for isocanted cubes.


52A38, Length, area, volume and convex sets

51A50, Polar geometry, symplectic spaces, orthogonal spaces

52A40, Inequalities and extremum problems involving convexity in convex geometry

\end{abstract}
\section{Introduction}\label{sec:intro}

Volume computation  is a meeting point of calculus and algebra, where integrals 
and determinants come into play. If through history,  only a few exact volume formulas of polytopes of higher dimension have been
proved, it
 is due to the fact that this question is, in general, very difficult. 
Since the 1980's it is known that exact volume computation  is
a $\#P$ problem; see the surveys \cite{Gritzmann_Klee_chap,Gritzmann_Klee,Khachiyan_chap} and the books \cite{Bisz_al,Matousek}. \label{citas_1}
A great number of approximated volume computation methods have been developed. The list of  references is enormous  and is out of the
scope of this paper.
Other approaches to volume computation are discrete volume, \emph{Brion theorem}  or  the exponential valuation; cf. \cite{Barvinok,Beck_Robins}.\label{citas_2} In \cite{Lawrence} Lawrence gives a remarkable exact volume formula  for simple $d$--polytopes. However, Lawrence formula is  not applicable here,  because the polytopes we study are not simple, but only almost--simple; { cf. Notation \ref{nota:common} for these notions}.

\bigskip

Some standard families of $d$--polytopes
are: simplices and  hypersimplices,  cubes and their polar duals, named cross--polytopes, cyclic polytopes, matroid polytopes and
order polytopes.
In  the papers \cite{Gritzmann_Klee,Kalai,Ziegler}\label{citas_3} the authors ask for more examples, i.e., more families of polytopes of
arbitrary dimension for which exact computations can
be performed and conjectures can be checked. We propose isocanted cubes for this matter.
 A number of concrete results (such as $f$--vector computation, semi--simplicity, cubicality  and the validity of
 several conjectures)  are proved in \cite{Puente_Claveria,Puente_Claveria_Iso}\label{citas_4} for isocanted cubes.\footnote{{In \cite{Puente_Claveria,Puente_Claveria_Iso}, \emph{isocanted cubes} are given the longer name of  \emph{isocanted alcoved polytopes}.}}


\bigskip

In $\R^d$ with $d\ge2$, what is a $d$--isocanted cube?
Given real parameters $0<a<\ell$, it is the convex set, denoted $\II_d(\ell,a)$,  defined by the following { inequalities}
\begin{equation}\label{eqn:isocanted_by_ineqs}
{ -\frac{\ell}{2}\le x_j\le \frac{\ell}{2}, \quad a-\ell\le x_j-x_k\le \ell-a,} \quad \forall j,k\in\{1,2,\ldots,d\}, j\neq k.
\end{equation}

{
For fixed dimension $d$, there exists only one combinatorial type of isocanted $d$--cubes, no matter the  values of $\ell$ and $a$.
This is proved in Corollary 3.13 in \cite{Puente_Claveria_Iso}.\label{citas_33} In fact, something stronger is proved there: the
face lattice of $\II_d(\ell,a)$ is equal to the lattice of proper subsets of $[d+1]$.
}

\begin{figure}[ht]
\centering
\includegraphics[width=6cm]{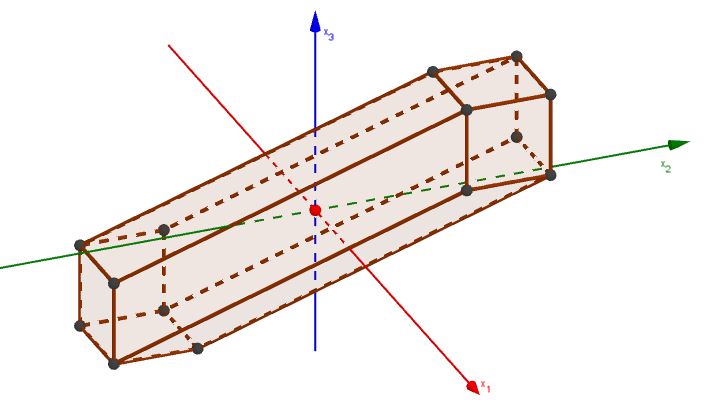}\ \includegraphics[width=6cm]{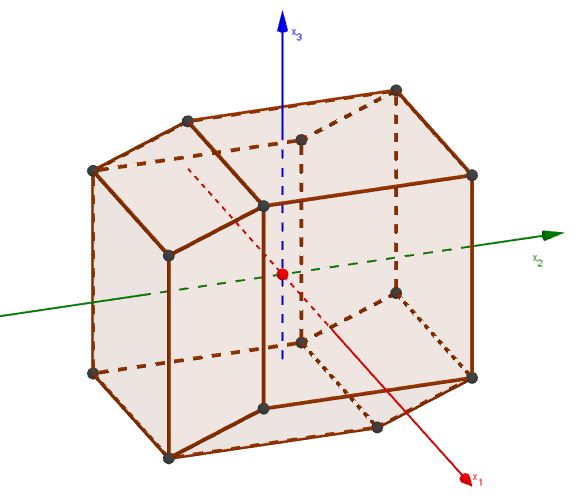}
\caption{Isocanted cube ($d=3$): cases { $\ell/2<a$} (left) and { $\ell/2>a$} (right).}\label{fig_isocanted_3D}
\end{figure}

To cant is a term taken from carpentry. It means to bevel, i.e., to change an edge of an object into a sloping face.
The case $d=3$ is shown in figure \ref{fig_isocanted_3D}: it is a polyhedron  combinatorially equivalent to a
\emph{rhombic dodecahedron}   (it has 14 vertices). 

\begin{figure}[ht]
\centering
\includegraphics[width=6cm]{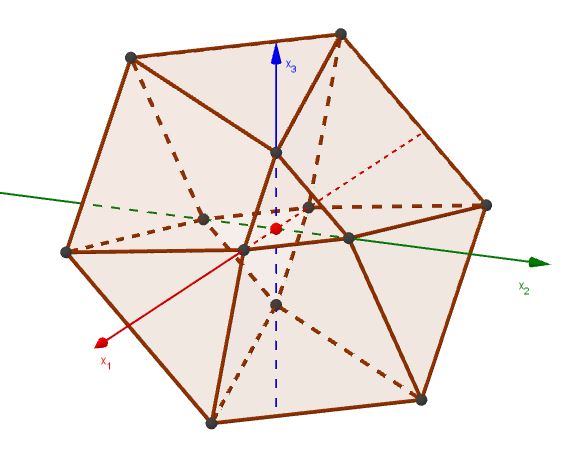}\ \includegraphics[width=6cm]{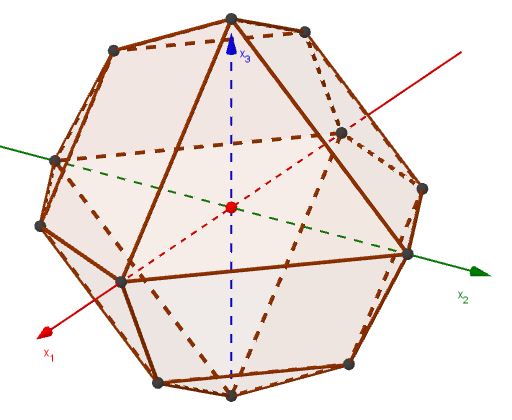}
\caption{Polar dual of the isocanted cube ($d=3$): cases { $\ell/2<a$} (left) and { $\ell/2>a$} (right).}\label{fig_polar_iso_3D}
\end{figure}

{The polar dual of the former, denoted $\II_d^\circ(\ell,a)$,  is  combinatorially equivalent
to a \emph{cuboctahedron}; it is shown in figure \ref{fig_polar_iso_3D}.}

\bigskip

Our main results are the two volume formulas presented in Theorem \ref{thm:vol_isobiselado} and Corollary \ref{cor:vol_polar_iso}.
The formulas are
\begin{equation}\label{eqn:v_I_d}
\vol_d\II_d(\ell,a)=(\ell-a)^{d-1}(\ell+(d-1)a)={\det\begin{pmatrix}
                                                      \ell & a & \cdots & \cdots & a \\
                                                      a & \ell & a & \cdots & a \\
                                                      \vdots & \ddots & \ddots & \ddots & \vdots \\
                                                      \vdots & &\ddots & \ddots  & a \\
                                                      a & \cdots & \cdots & a & \ell
                                                    \end{pmatrix}}
\end{equation}
and
\begin{equation}\label{eqn:v_I_d_o}
\vol_d \II_d^\circ(\ell,a)= \frac{2^{d+1}}{\ell^d d!}\sum_{j=0}^{d-1}{{d+j-1}\choose{j}}\left(\frac{\ell}{2(\ell-a)}\right)^j.
\end{equation}
The matrix in (\ref{eqn:v_I_d}) was introduced by Bose, (cf.\cite{Babai,Bose}).\label{citas_Bose} Formulas (\ref{eqn:v_I_d}) and (\ref{eqn:v_I_d_o}) specialize to $a=0$ properly.  Indeed, $\II_d(\ell,0)=C_d(\ell)$ is  the $d$--cube of edge--length $\ell$, given by facet equations $x_i=\pm\frac{\ell}{2}, i=1,2,\ldots,d$. Its polar dual is $C_d^\circ(\ell)=\conv(\pm\frac{2}{\ell}e_i:i=1,2,\ldots,d)$,  a $d$--cross--polytope of edge--length equal to $\frac{2\sqrt2}{\ell}$. \label{polar_of_cube} It is known { (cf. the second line in Table 15.2.1 in \cite{Handbook_chap15}, where the volume of $X_d(\sqrt2)$ is given as $\frac{2^d}{d!}$)} that the volume of a $d$--cross polytope of edge--length $\ell'$, denoted $X_d(\ell')$, is given by the formula
\begin{equation}\label{eqn:X}
 \vol_d X_d(\ell')= \frac{1}{d!}\sqrt{2^d}{\ell'}^d.
\end{equation}
Substituting $\ell'=\frac{2\sqrt2}{\ell}$ in (\ref{eqn:X}) and  using {the so called \emph{finite Companion Binomial
Theorem}} (cf. (\ref{eqn:finite_companion_binomial}) below)
we get
\begin{equation}
 \vol_d X_d
 \left(\frac{2\sqrt2}{\ell}\right)= \frac{4^d}{\ell^dd!}=\frac{2^{d+1}}{\ell^d d!}\sum_{j=0}^{d-1}{{d+j-1}\choose{j}}\left(\frac{1}{2}\right)^j=\vol_d \II_d^\circ (\ell,0).
\end{equation}

A new notion  in this paper is roof,
introduced in Definition \ref{dfn:right_roof}. Euclid used certain roofs in his construction of
a regular dodecahedron from a cube (cf. \cite{Cromwell, Sommerville}).\label{citas_Euclid} Roofs form a subclass of prismatoids. Roofs are  important for us because
we decompose $\II_d^\circ(\ell,a)$ into a finite union of pyramids whose bases are roofs. The volume of a $d$--pyramid based on
a $(d-1)$--{ polytope} $P$ (with apex $Q$ and height $h$) is known to be
\begin{equation}\label{eqn:vol_pyramid_intro}
    \vol_d \pyr(P,Q)=\frac{(\vol_{d-1} P) h}d.
    \end{equation}
  { The definition of roof is a subtle one, because it embodies combinatorics as well as metric  facts.} In  Lemma \ref{lem:vol_roof} we compute the volume of a roof
by  elementary methods (integration of hyperplane sections).  Particular cases of $d$--roofs are: isosceles triangles and  isosceles
trapezia ($d=2$) and
right pyramids (equilateral--triangular--based or square--based), frusta of right  equilateral--triangular--based pyramids,
right prisms (equilateral--triangular--based or square--based) ($d=3$).\footnote{ Properly, we should speak of \emph{right roofs.} We write just roofs for simplicity.}
Roofs have one very special facet, called \emph{major { base}}. The major { base} is the product of two regular simplices. This is why the
volume formula for a roof
is derived from the known volume formula for a regular simplex (of edge--length $\ell$, denoted $\Delta_d(\ell)$) and the product formula
\begin{equation}\label{eqn:vol_simplex_and_prod}
\vol_d \Delta_d(\ell)=\frac{1}{d!}\sqrt{\frac{d+1}{2^d}}\ell^d,\qquad  \vol (P \times P')=\vol (P) \vol (P'),
\end{equation}
{ (cf. the third line in Table 15.2.1 in \cite{Handbook_chap15}, where the volume of $\Delta_d(\sqrt2)$ is given as $\frac{\sqrt{d+1}}{d!}$).}
\bigskip

In the centrally symmetric case, \emph{the Mahler conjecture} asks whether the  product of the volumes of a  centrally symmetric
$d$--convex body and of its  polar dual is bounded
below by the product of volumes of a $d$--cube and of its polar dual.
For $C_d(\ell)$ and its polar dual { $X_d\left(\frac{2\sqrt2}{\ell}\right)$},  the volume product is equal to
\begin{equation}\label{eqn:prod_vol_cube}
{\vol_d C_d(\ell) \vol_d X_d\left(\frac{2\sqrt2}{\ell}\right) }=\ell^d \frac{4^d}{\ell^d d!}=\frac{4^d}{d!}
\end{equation}
and does not depend on $\ell$.
The conjecture has
been proved
in the case $d=2$ by Mahler in \cite{Mahler}, in the case  $d=3$ by Iriyeh and Shibata
in \cite{Iriyeh_Shibata} (a shorter proof is found in \cite{Fradelizi_al}), and in arbitrary dimension for some  classes of convex
bodies, such as zonoids, in \cite{Gordon_al,Reisner}, other particular
cases \cite{Barthe_Fradelizi, Lopez_Reisner}. \label{citas_5} It has been proved
in \cite{Reisner_al} that bodies having a boundary point where the generalized Gauss curvature is positive are not  local minimizers.
In \cite{Kim} it has been shown that  \emph{Hanner polytopes} are local   minimizers, and in \cite{Nazarov_al} \label{citas_6} that the unit cube is
a strict local minimizer with respect to the \emph{Banach--Mazur distance}. In relation with  normed spaces, the Mahler conjecture has
been confirmed for
the unit ball of a finite real \emph{Banach space} with a 1--unconditional basis, in  \cite{Meyer,Reisner_uncondit}. For the unit ball
of  the Lipschitz--free space arising from a  finite pointed  metric space, partial results about the conjecture have been obtained in \cite{Alexander_al,Godard}. In \cite{Karasev} Mahler conjecture has been proved for some hyperplane sections of $\ell_p$ balls
(and projections), $1\le p\le +\infty$,  and of Hanner polytopes. \label{citas_7}
 The  interest on this conjecture and related inequalities remains uninterrupted, and here we have cited only a few references.
 See bibliography in these papers, as well as the theses \cite{Henze,Hupp} for a more complete view.\label{citas_8}

\bigskip

$\II_d(\ell,a)$ is the sum of $d+1$ segments, whence it is a zonoid. Thus, the Mahler conjecture holds true for $\II_d(\ell,a)$.
The polar dual of a zonoid is a  central section of the unit ball of a finite $\ell_1$--space (cf. \cite{Bolker}).
The Mahler conjecture for sections as such, when arising from finite metric spaces, has been recently studied in several papers.
In Remarks \ref{rem:M} and \ref{rem:M_2},  we have explained the relation between $\II_d(\ell,a)$ and  some results in  \cite{Alexander_al,Godard}.\label{citas_9}

\bigskip

The paper is organized as follows. In Section \ref{sec:ring_of_mat} we study the matrices of the form $\alpha I_d+\beta J_d$.
In Section \ref{sec:vol_isocanted}, we prove formula (\ref{eqn:v_I_d}), using Section \ref{sec:ring_of_mat}.
Section \ref{sec:background} gathers some well--known facts.
{In Section \ref{sec:ring_of_mat_revisited} a relationship between polar duals and inverse matrices  appears as a  by--product: to  the matrix
${M_d(\ell,a)}=(\ell-a)I_d+aJ_d$ we associate,
in a unique way,   an origin--symmetric parallelepiped, denoted $\tau (Par({M_d(\ell,a)}))\subseteq \R^d$. Then, the vertices of its polar dual
are the columns of $\pm 2({M_d(\ell,a)})^{-1}$.}  In
Section \ref{sec:vol_polar_isocanted} we compute the volume of a roof and prove formula (\ref{eqn:v_I_d_o}). Finally, in
Section \ref{sec:Mahler} we reformulate
\emph{the Mahler conjecture} (for $\II_d(\ell,a)$) as the question of positivity on the
open interval $(0,1)$ of  a polynomial $p_d(x)$ of degree $d$ in one variable. We prove $p_d(x)>0$ only using the known bound for
central binomial coefficients (cf. \cite{Wikipedia_central_bin})\label{citas_10}
\begin{equation}\label{eqn:bound_central}
{{2d-2}\choose {d-1}} \ge\frac{4^{d-1}}{\sqrt{\pi\left(d-\frac{1}{2}\right)}}
\end{equation}
and the \emph{Descartes rule of signs}.

\bigskip

Alcoved polytopes of affine type A
arising from affine Coxeter arrangements are studied in \cite{Lam_Postnikov} and a subsequent paper, and more recently in \cite{Sjoberg}.
These are lattice polytopes
over $\Z^d$. In \cite{Werner_Yu}, \label{citas_11} symmetric alcoved polytopes of the
same kind are studied in relation with tropical geometry. A combinatorial volume formula is given in \cite{Lam_Postnikov}. {Isocanted cubes are a two parameter example of alcoved polytopes, as defined in \cite{Puente_QE}.\label{citas_32}} Despite the same name, this volume formula is not applicable to
isocanted cubes, since isocanted cubes are not,
in general, lattice polytopes.

\bigskip

The calculation of volumes is present in Mathematics from ancient times. People were interested in computing the content of receptacles,
the amount of bricks needed for a construction, etc.  The oldest known volume formula is Egyptian,  found in Problem 14 in the
\emph{Moscow Mathematical Papyrus} (ca. 1850 BC). It  expresses the
volume of a frustum of a  pyramid   as
\begin{equation}\label{eqn:egypt}
  h\frac{a^2+ab+b^2}{3},
\end{equation} where
 $a^2$, $b^2$ are the areas of the frustum bases and $h$ is  the height (cf. \cite{Boyer}\label{citas_30}). {From the
 papyrus text and figure, it is unclear whether the pyramid is right or not as well as what is the shape of the bases.}
 We recover this formula (when the bases are equilateral triangles { and the pyramid is right}), as  a particular  case of
 roof (cf. (\ref{eqn:egypt_2})).

\bigskip

Notations: $d\ge2$ denotes the dimension of the ambient space, $O$ denotes the origin in $\R^d$,
$(e_1,e_2,\ldots,e_d)$ denotes the standard { vector} basis. Capital letters  $A,B,D,H,I,J,M,S$ denote matrices,
$K\subseteq \R^d$ denotes a convex body, $P,Q$  denote polytopes,  $H, H_1,H_2, MH, mh$ denote hyperplanes,
$\ell,\ell_1,\ell_2,\ell_3,a,b,c,h$ denote positive real numbers,
$p,r$ denote polynomials (in one or two variables), Greek letters
$\alpha, \beta,\gamma, \delta$ denote elements in  ${\bf k}$, a subring of $\R$,
$\phi, \rho,\sigma, \tau:\R^d \rightarrow \R^d$ denote
maps, $V,C$ denote positive integers, $v,v'$ denote vertices of a polytope, { $x,y,z,p_1,p_2,\ldots$ denote points in $\R^d$. We identify $d$--tuples with coordinates of points  or vectors in $\R^d$, in the usual way.} The set $\{1,2,\ldots,d\}$ is denoted $[d]$.

\section{Two rings of matrices and  many equal--volumed parallelepipeds}\label{sec:ring_of_mat}
Let ${\bf k}$ be a subring of $\R$. 
For $d\in\N$, let   $\MM_d^{sc}({\bf k})$ be the  commutative ring of $d\times d$
symmetric circulant  matrices over ${\bf k}$.\footnote{The product of circulant matrices is circulant. Any two circulant matrices commute. The product of two commuting symmetric matrices is symmetric. For information on circulant matrices, cf. \cite{Davis,Gray}.\label{citas_12}} Let $I_d$ and $J_d$ be the  identity matrix and the all--ones matrix, respectively.
Since  $J_d^2=dJ_d$, then the set
\begin{equation}
\left\{\alpha I_d+\beta J_d: \alpha, \beta\in {\bf k}\right\}
\end{equation}
is easily checked to be a subring of  $\MM_d^{sc}({\bf k})$. The product formula is
\begin{equation}\label{eqn:product_formula}
\left(\alpha I_d+\beta J_d\right)\left(\gamma I_d+\delta J_d\right)=\alpha\gamma I_d+\left(\alpha\delta+\beta\gamma+d\beta\delta\right) J_d.
\end{equation}
The 
spectrum of $\alpha I_d+\beta J_d$ is
\begin{equation}
\left\{\alpha \text{\ with\ mult.\ }(d-1), \alpha+d\beta\right\}.
\end{equation}
If $\beta\neq0$, then the matrix $\alpha I_d+\beta J_d$ is diagonalizable over $\R$ (because it is symmetric and $\dim \ker (\alpha I_d+\beta J_d-\alpha I_d)=\dim \ker J_d=d-1$). We get
\begin{equation}\label{eqn:diagonalization}
\alpha I_d+\beta J_d=H_d D_d(\alpha,\beta) H_d^{-1}
\end{equation}
with
\begin{equation}\label{eqn:diagonal}
D_d(\alpha,\beta):=\diag\left(\alpha,\ldots,\alpha,\alpha+d\beta\right)
\end{equation}
and the so called \emph{Helmert matrix} (cf. p.5 in \cite{Lancaster} \label{citas_13})
\begin{equation*}\label{eqn:Helmert}
H_d=
\begin{pmatrix}
1/\lambda_{1}&1/\lambda_{2}&1/\lambda_{3}&\cdots&1/\lambda_{d-1}&1/\sqrt{d}\\
-1/\lambda_{1}&1/\lambda_{2}&1/\lambda_{3}&\cdots&1/\lambda_{d-1}&1/\sqrt{d}\\
0&-2/\lambda_{2}&1/\lambda_{3}&\cdots&1/\lambda_{d-1}&1/\sqrt{d}\\
0&0&-3/\lambda_{3}&&\vdots&\vdots\\
\vdots&\vdots&&&1/\lambda_{d-1}&1/\sqrt{d}\\
0&0&0&&-(d-1){/}\lambda_{d-1}&1/\sqrt{d}\\
\end{pmatrix}
\end{equation*}
and $\lambda_j:=\sqrt{j+j^2}$. It is easy to check that  $H_d$ is orthogonal (i.e., $H_d^T=H_d^{-1}$) and has zeros below the
first subdiagonal.

If $\alpha(\alpha+d\beta)\neq0$, then  $\alpha I_d+\beta J_d$ is invertible. { Indeed,}
\begin{equation}
{\left(D_d(\alpha,\beta)\right)^{-1}=\diag(1/\alpha,\ldots,1/\alpha, 1/(\alpha+d\beta))=D_d(1/\alpha,\beta')}
\end{equation}
\begin{equation}\label{eqn:inverse}
\left(\alpha I_d+\beta J_d\right)^{-1}=H_d \left(D_d(\alpha,\beta)\right)^{-1} H_d^{-1}=H_d D_d(1/\alpha,\beta') H_d^{-1}=\frac{1}{\alpha}I_d+\beta' J_d,
\end{equation}
with  {  $\frac{1}{\alpha}+d\beta'=\frac{1}{\alpha+d\beta}$ and  $\beta'=\frac{-\beta}{\alpha(\alpha+d\beta)}$ where the equality  $\left(\alpha I_d+\beta J_d\right)\left(\frac{1}{\alpha}I_d+\beta' J_d\right)=I_d$ is readily checked using (\ref{eqn:product_formula}).}

{
\begin{nota}[Parallelepiped determined by a square matrix]\label{nota:par}
The vector columns of a $d\times d$   matrix $M$ over ${\bf k}$ determine the  edges of the parallelepiped $Par(M)\subseteq \R^d$ having one vertex at the origin $O$.
\end{nota}
}

\begin{rem}[Equal volumed parallelepipeds]\label{rem:parallelepipeds}
 The
 $d$--volume of $Par(M)$ is $|\det M|$. It follows from (\ref{eqn:diagonalization}) that the parallelepipeds
 $Par(\alpha I_d+\beta J_d)$, $Par(D_d(\alpha,\beta))$ and $Par(SD_d(\alpha,\beta)S^{-1})$ all
 have the same $d$--volume, for each  regular matrix $S\in\MM_d({\bf k})$.
\end{rem}

\begin{rem}[Rotation--invariant parallelepiped]\label{rem:parallelepipeds_2}
Consider the rotation
\begin{equation}\label{eqn:rotation}
\sigma_{2\pi/d}:\R^d\rightarrow\R^d
\end{equation} around the axis $x_1=x_2=\cdots=x_d$ with angle $2\pi/d$.
The
circulant matrix ${B_d}$ (having first column equal to $(0,1,0,\ldots,0)^T$)  represents $\sigma_{2\pi/d}$ with respect
to the standard  { vector} basis $(e_1,e_2,\ldots,e_d)$.
The set of columns of { the product $B_d (\alpha I_d+\beta J_d)$} is equal to the set of columns of $\alpha I_d+\beta J_d$, meaning that the parallelepiped $Par(\alpha I_d+\beta J_d)$ is invariant under $\sigma_{2\pi/d}$.
\end{rem}


\section{Volume of the isocanted cube}\label{sec:vol_isocanted}

In $\R^d$ consider the $d$--cube of
edge--length $\ell>0$, centered at the origin $O$ and having edges parallel
to the coordinate axes. For 
$0<a<\ell$, consider the $d$--polytope  $\II_d(\ell,a)$ {given by the  linear inequalities in (\ref{eqn:isocanted_by_ineqs}).}

\begin{figure}[ht]
\centering
\includegraphics[width=6cm]{fig_isocanted_2D_big_corner_v2.jpg}\ \includegraphics[width=6cm]{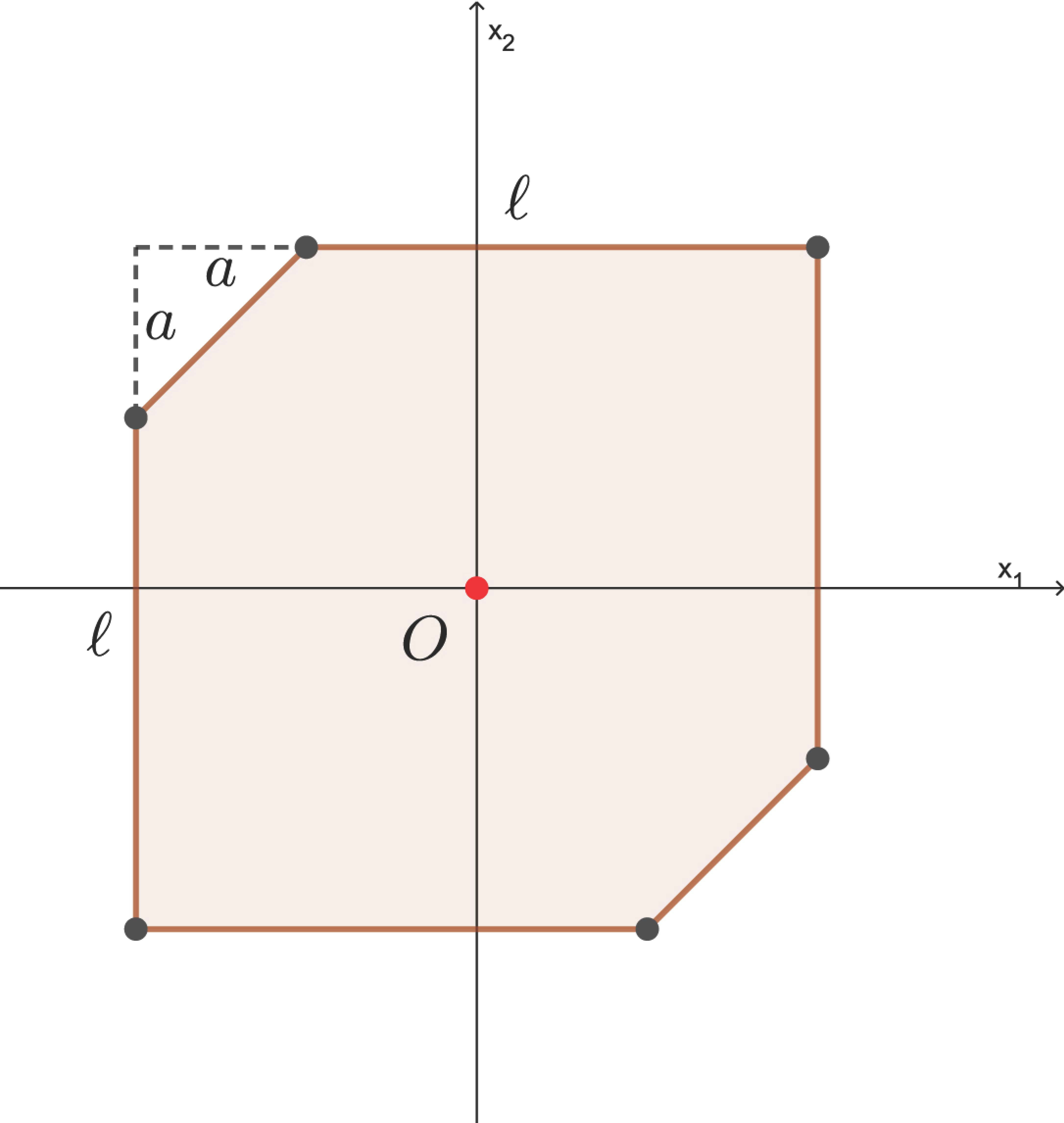}
\caption{Isocanted {square} ($d=2$): cases { $\ell/2<a$} (left) and { $\ell/2>a$} (right).}\label{fig_isocanted_2D}
\end{figure}

We say that $\II_d(\ell,a)$ is the \emph{isocanted cube} 
with \emph{cant parameter} $a$. 
See figures \ref{fig_isocanted_2D} and  \ref{fig_isocanted_3D}. 
In \cite{Puente_Claveria_Iso}\label{citas_14} we proved that isocanted cubes are
zonotopes.  Moreover, in the proof of Theorem \ref{thm:vol_isobiselado} below, we give an explicit  sum description.qual
Clearly, $\II_d(\ell,a)$ is symmetric with respect to the origin $O$.
We write $\II_d$ whenever $\ell,a$ are understood. Let

\begin{equation}\label{eqn:M}
M_d(\ell,a):=(\ell -a)I_d+aJ_d=\begin{pmatrix}
                                                      \ell & a & \cdots & \cdots & a \\
                                                      a & \ell & a & \cdots & a \\
                                                      \vdots & \ddots & \ddots & \ddots & \vdots \\
                                                      \vdots & &\ddots & \ddots  & a \\
                                                      a & \cdots & \cdots & a & \ell
                                                    \end{pmatrix}
\end{equation}

\begin{thm}[Volume of the isocanted cube]\label{thm:vol_isobiselado}
For each dimension $d\ge2$ and parameters $0<a<\ell$, the $d$--volume of $\II_d(\ell,a)$ is equal to the determinant of the \emph{Bose matrix}
\begin{equation}\label{eqn:r}
r_d(\ell,a):=(\ell-a)^{d-1}(\ell+(d-1)a)=(\ell-a)^{d}+da(\ell-a)^{d-1}=\det M_d(\ell,a).
\end{equation}

\end{thm}
\begin{proof}
The second equality in (\ref{eqn:r}) is straightforward and the last one follows from Section \ref{sec:ring_of_mat}.
The polytope $\II_d(\ell,a)$ is the \emph{Minkowski sum} (i.e.,  vector sum)  of $d+1$ segments,
namely
\begin{equation}\label{eqn:zonotope}
\II_d(\ell,a)=\sum_{i=1}^{d+1}\left[-y_i,y_i\right],\qquad \text{with \ } y_i=\begin{cases}
                                                                                (\ell-a)e_i/2, & \mbox{if } i\in [d] \\
                                                                                {a}(e_1+e_2+\cdots+e_d)/2,
                                                                                & \mbox{if }i=d+1.
                                                                                \end{cases}
\end{equation}
{  The equality (\ref{eqn:zonotope}) holds indeed: a point $z$ in $\sum_{i=1}^{d+1}\left[0,y_i\right]$
satisfies $z=\left(\sum_{i=1}^{d}\mu_i e_i\right)+\mu_{d+1}(e_1+e_2+\cdots+e_d)=\sum_{i=1}^{d}(\mu_i+\mu_{d+1}) e_i$,
with $0\le\mu_i\le \frac{\ell-a}{2}$ and $0\le\mu_{d+1}\le \frac{a}{2}$, and so the coordinates of $z$  satisfy
inequalities (\ref{eqn:isocanted_by_ineqs}).
By the symmetry of the set on the left--hand--side of (\ref{eqn:zonotope}), this is enough to prove that the left--hand--side is contained in the right--hand--side. The converse containment is proved similarly.
}

Then, the \emph{McMullen--Shephard} formula (cf. (58) in \cite{Shephard}\label{citas_15} and also \cite{Montgomery_al} and Corollary 3.4 in \cite{Gover_Krikorian}) ---providing the exact volume for a zonotope--- applies, and we get
\begin{equation}\label{eqn:shephard}
\vol_d \II_d(\ell,a)={2^d}\sum_{i=1}^{d+1} |\det M_i|
\end{equation}
where $|\det M_i|$ is the absolute value of the $d\times d$ minor obtained from the following   $d \times (d+1)$ matrix by omitting  the $i$--th column
\begin{equation*}
\left(\begin{array}{ccccc}
\frac{\ell-a}{2}&0&\cdots&0&\frac{a}{2}\\
0&\frac{\ell-a}{2}&&0&\frac{a}{2}\\
\vdots&&\ddots&\ddots&\vdots\\
0&\vdots&\frac{\ell-a}{2}&0&\frac{a}{2}\\
0&\cdots&0&\frac{\ell-a}{2}&\frac{a}{2}\\
\end{array}\right).
\end{equation*}
If we omit the last  column, we get a  diagonal matrix and $|\det M_{d+1}|=\left(\frac{\ell-a}{2}\right)^d$ and the remaining  summands $|\det M_i|$ are  {also diagonal, up to column permutation, and are} all equal to  $\frac{a\left(\ell-a\right)^{d-1}}{2^d}$.  {The left--hand--side of (\ref{eqn:shephard}) ammounts to $(\ell-a)^{d}+da(\ell-a)^{d-1}=r_d(\ell,a)$, as} desired.
\end{proof}

We also have a proof of the former result by integration and induction on the dimension.

{
Application:
$d\ge2$ friends have an appointment between 10h and 11h. Each one arrives at random in the given time interval and waits for only 10 minutes.  What is the probability that they all meet together?

Solution: The total sample space is the unit $d$--cube $C_d(1)$.
The event is $\II_d(1,a)$ with
$1-a=\frac{10}{60}$ so $a=\frac{5}{6}$ and
 the  probability  is
\begin{equation}\label{eqn:probability}
\vol_d\II_d\left(1,\frac{5}{6}\right)=\left(1-\frac{5}{6}\right)^{d-1}\!\!\left(1+(d-1)\frac{5}{6}\right)=\!\frac{5d+1}{6^d}.
\end{equation}
Thus, if $d$ is large, it is very rare that all $d$ meet together!
}

\begin{rem}\label{rem:properties_r}
$r_d(\ell,a)$ in (\ref{eqn:r}) is an instance of  { a} homogeneous polynomial of  degree $d$  in $\ell,a$ with coefficients
in ${\bf k}=\Z$ which is equal to the determinant of a symmetric matrix over { the ring} ${\bf k}[\ell,a]$ with
homogeneous linear entries.

The general question of \emph{determinantal hypersurfaces} (i.e., which homogeneous polynomials are the determinant of a matrix
(possibly symmetric) with homogeneous entries) is classical in Algebraic Geometry. A modern approach (including historical comments
and references) is found in \cite{Beauville}.\label{citas_16} It studies generic forms over ${\bf k}$, an arbitrary field, not
necessarily algebraically closed. Notice that  $r_d(\ell,a)$  in (\ref{eqn:r}) is the product of $d$ linear factors over $\Z$,
and so it is very far from being generic.
\end{rem}

\begin{question}
Given $d\ge2$, can we find more families of $d$--polytopes and of $d\times d$ matrices, depending on some real parameters ($\ell,a$ etc.) such that the  $d$--volume of one such polytope is the determinant of one such matrix?
\end{question}

\section{Background}\label{sec:background}
Recall several commonly known facts on polytopes,
in the affine and affine--Euclidean setting $\R^d$; cf. \cite{Blind_almost_simple,Coxeter,Handbook_chap15,Sommerville,Ziegler_book}.\label{citas_17}

\begin{nota}[Common polytopes and common operations]\label{nota:common}
\begin{itemize}
\item  $P\times P'$ denotes the \emph{Cartesian product}  of the polytopes $P$, $P'$.
  If $P\subseteq \R^d$ and $P'\subseteq \R^{d'}$, then the coordinates of the vertices of $P\times P'\subseteq \R^{d+d'}$ are obtained by concatenation of the $d$ coordinates { of} a vertex in $P$ with the $d'$ coordinates   { of} a vertex in $P'$ in all possible ways. \label{item:producto}  It holds
\begin{equation}\label{eqn:vol_product}
\dim (P\times P')=\dim (P)+  \dim (P'), \qquad \vol (P \times P')=\vol (P) \vol (P').
\end{equation}

\item Given a $(d-1)$--polytope $P$ and a point $Q\not\in\aff(P)$, let  $\pyr(P,Q)$ be  the \emph{pyramid} with  \emph{{ base}} $P$
and \emph{apex} $Q$.  We have
    \begin{equation}\label{eqn:vol_pyramid}
    \vol_d \pyr(P,Q)=\frac{(\vol_{d-1} P)  \dd(\aff(P),Q)}d.
    \end{equation}
    and $\dim \pyr(P,Q)=1+\dim P$.
{
\item   $\conv(x,y):=\{\lambda x+(1-\lambda)y\in\R^d: 0\le \lambda\le 1\}$ denotes the \emph{(line) segment} determined by
points $x,y\in \R^d$. Another notation is $[x,y]$.
\item A $d$--polytope is \emph{regular} if its facets  and its half--edge vertex figures are all regular.

\item The term \emph{right} applied to a polytope means orthogonality of some of its elements. For instance, a pyramid is right if
the line segment joining its apex with the barycenter of its  base is orthogonal to the base.

\item  The term \emph{isosceles} applied to a polytope means equal length of some of its elements. For instance, a pyramid is isosceles if
the line segments joining its apex with each vertex in the base have all equal length.
}
\item  ${\frac{\ell}{\sqrt2}\conv(e_1,e_2,\ldots,e_d)}$ is  the   \emph{standard regular simplex} of dimension $d-1$ and edge--length $\ell\ge0$. { It is contained in the affine hyperplane of equation $x_1+x_2+\cdots+x_d=\frac{\ell}{\sqrt2}$.}
\item By \emph{isosceles}   $d$--simplex we mean  a pyramid whose { base} is a regular $(d-1)$--simplex  and whose remaining edges have all the same length.



{

\item A \emph{parallelogram} is obtained by translating a segment by a vector  not parallel to the segment. A $d$--\emph{parallelepiped} is the result of  the analogous $d$--dimensional construction. \footnote{ Some authors (not all) use the term \emph{parallelotope} as a  synonym.}
\item A $d$--\emph{box} (or $d$--\emph{orthotope}) is  a $d$--parallelepiped having mutually orthogonal facets.
\item A $d$--\emph{cube} is a  $d$--box all whose edges have equal length.


\item A $d$--polytope is \emph{simple} if exactly $d$ facets meet at each of its vertices.
\item A $d$--polytope is \emph{almost simple} if exactly $d$ or $d+1$ facets meet at each of its vertices.
\item The polar dual of a simple polytope is a \emph{simplicial} polytope.
\item The \emph{valency} (or \emph{degree}) of a vertex is the number of facets containing it.


\item A polytope is \emph{cubical} if each  proper face in it is a cube.
}
\end{itemize}
\end{nota}

\begin{nota}[$\Delta_{d-1}(\ell)$]\label{nota:Delta}
{For the rest of the paper, we fix one rigid motion in $\R^d$ such that the image of the simplex $\frac{\ell}{\sqrt2}\conv(e_1,e_2,\ldots,e_d)$ is contained in the hyperplane of equation $x_d=0$ and centered at the origin. This image of $\frac{\ell}{\sqrt2}\conv(e_1,e_2,\ldots,e_d)$ is denoted
$\Delta_{d-1}(\ell)$. By abuse of notation, the translation of $\Delta_{d-1}(\ell)$ to any parallel  hyperplane  and centered at the origin of that hyperplane is also denoted $\Delta_{d-1}(\ell)$.
}

\end{nota}
\begin{nota}[Congruence]\label{nota:congruent}
Two polytopes $P,Q\subseteq \R^d$ are \emph{congruent} if there exists a rigid motion $\sigma:\R^d\rightarrow \R^d$ such that $\sigma(P)=Q$.
We will write $P\simeq Q$.\footnote{ A rigid motion is an isometry and therefore congruent polytopes are isometric.}

\end{nota}

Recall  three well--known volume formulas: $d$--cube, $d$--cross polytope and $d$--regular simplex with edge--length equal to $\ell$:
\begin{equation}\label{eqn:known}
\ell^d,\qquad \frac{1}{d!}\sqrt{2^d}\ell^d,\qquad\frac{1}{d!}\sqrt{\frac{d+1}{2^d}}\ell^d.
\end{equation}
\medskip

For each $j,k\in\N$, and each $h\ge0$, the following integral\footnote{This integral is related to the \emph{binomial distribution}.}
is proportional to $h$ with a factor given by
 the \emph{beta function} (cf. p.131 in \cite{Sommerville}\label{citas_18}
 (3.8) in \cite{Gould_III}, or p.215 in \cite{Webster}):
\begin{equation}\label{eqn:integrales}
\int_{0}^h\left(1-\frac{s}{h} \right)^j\left(\frac{s}{h} \right)^kds=h \beta(j+1,k+1),\ \text{where\ } \sum_{n=0}^{j}{{j}\choose{n}}\frac{(-1)^n}{k+n+1}=\beta(j+1,k+1)
={\frac{j!k!}{(j+k+1)!}}\in \Q.
\end{equation}

\bigskip

The \emph{Lagrange  distance  formula} between a point $P\in \R^d$  with coordinates  $(p_1,\ldots, p_d)$ and a hyperplane $H\subseteq \R^d$ with equation
$a_1x_1+\cdots {+} a_dx_d+a_0=0$ yields
\begin{equation}\label{eqn:lagrange}
\dd(P,H)=|a_1p_1+\cdots {+} a_dp_d+a_0|/\sqrt{a_1^2 +\cdots+a_d^2}.
\end{equation}

\bigskip

The \emph{finite Companion Binomial} is
(cf. formula (1.16) in \cite{Gould_II}\label{citas_31})
\begin{equation}\label{eqn:finite_companion_binomial}
  \sum_{j=0}^{d-1}{{d+j-1}\choose{j}}\left(\frac{1}{2}\right)^j=2^{d-1}.
\end{equation}

\bigskip

For $d,j,r$ non--negative integers with $0\le j\le d+r$ the following \emph{combinatorial equality} holds
\begin{equation}\label{eqn:igualdad_combi}
{{d+r}\choose{j} }=\sum_{0\le n\le r; 0\le j-n\le d}{{d}\choose{j-n} }{{r}\choose{n} }=\sum_{n=0}^{\min\{r,j\}}{{d}\choose{j-n} }{{r}\choose{n} }
\end{equation}
the proof of which comes from the
observation  that  choosing $j$ elements from a set  of  $d+r$  given elements is achieved by
choosing $n$ elements from $r$ given ones,  and choosing $j-n$ from $d$ given ones,  for all possible values of  $n$.


{
A  combinatorial equality that we will use in Remark \ref{rem:Alexander} is
\begin{equation}\label{eqn:mi_combinatorial_eq}
2 \sum _{j=0}^{d-1}  {{d+j-1 }\choose j} =  {{2d}\choose d}.
\end{equation}
A proof of (\ref{eqn:mi_combinatorial_eq}) is given here: first, the so called \emph{Hockey stick pattern } ${{m}\choose{n}}={{m-n-1}\choose{0}}+{{m-n}\choose{1}}+\cdots+{{m-2}\choose{n-1}}+{{m-1}\choose{n}}=\sum_{j=0}^{n}{{m-n+j-1}\choose{j}}$ for all $n\le m-1$, is obtained by repeated use of the basic equality of combinatorial numbers: ${{m}\choose{n}}={{m-1}\choose{n-1}}+{{m-1}\choose{n}}={{m-2}\choose{n-2}}+{{m-2}\choose{n-1}}+{{m-1}\choose{n}} =\cdots$. Then we take
$m=2d-1$ and $n=d-1$, getting ${{2d-1}\choose{d-1}}=\sum_{j=0}^{d-1}{{d+j-1}\choose{j}}$, and then duplicate, getting ${{2d}\choose d}={{2d-1}\choose {d-1}}+{{2d-1}\choose {d}}=2{{2d-1}\choose {d-1}}=2\sum _{j=0}^{d-1}  {{d+j-1 }\choose j} $.

}

\bigskip
The \emph{polar dual} of  a set $K\subseteq \R^d$ is defined in terms of the usual inner product $\langle , \rangle$ as follows:
\begin{equation}\label{eqn:polar}
K^\circ:=\{x\in \R^d:\langle x,y\rangle\le 1, \forall y\in K\}
\end{equation}
In particular, when $K{=\conv(p_1,p_2,\ldots,p_r)}$ is the convex hull of a  finite set of points
 $p_1,p_2,\ldots,p_r\in \R^d$, then
\begin{equation}\label{eqn:polar_2}
K^\circ=\{x\in \R^d:\langle x,p_j\rangle\le 1, \forall j\in[r]\}.
\end{equation}
It is well--known that
\begin{equation}\label{eqn:bi_polar}
(K^\circ)^\circ=K,
\end{equation} whenever $K$ is closed, convex and contains the origin $O$ in its interior. { In particular, with the usual identification of points  and vectors in $\R^d$,  if
\begin{equation}\label{eqn:polar_3}
K=\{x\in \R^d:\langle x,u_j\rangle\le 1, \forall j\in[r]\}\qquad \text{ \ then \ }K^\circ=\conv(p_1,p_2,\ldots,p_r).
\end{equation}
}
See  \cite{Barvinok,Grunbaum,Handbook_chap15,Webster}\label{citas_19} for further details on polar duals.

\bigskip

{ 
A \emph{convex body} in  Euclidean space is a compact convex subset having non--empty interior. Consider a convex body $K\subseteq \R^d$ containing the origin $O$ in its interior.  The \emph{volume product} (also called \emph{Mahler volume})  of  $K$ is, by definition,
\begin{equation}\label{eqn:prod_vol}
\PP(K):=\vol_d (K)\vol_d (K^\circ)=\PP(K^\circ)
\end{equation}
Clearly, $\PP(K)=\PP(\sigma(K))$, for each invertible  linear map $\sigma:\R^d \rightarrow\R^d$.

\bigskip
If it were true that for given convex $d$--bodies $K_1,K_2$ the equality  $\vol_d(K_1)=\vol_d(K_2)$
implies $\vol_d(K_1^\circ)=\vol_d(K_2^\circ)$, then we would no be studying \emph{the Mahler conjecture}.
However,  this implication is false: take $K_1$ be the Euclidean unit ball and $K_2$ a $d$--cube of edge--length
$\ell>0$  chosen in  such a way that $\vol_d(K_1)=\vol_d(K_2)$. Then $K_1^\circ=K_1$ but (\ref{eqn:known}) shows that  $\vol_d(K_2)\neq \vol_d(K_2^\circ )$.

\begin{dfn}[Hanner polytope]\label{dfn:Hanner}
A polytope is \emph{Hanner} if it is obtained from segments, by a finite sequence of two transformations: taking Cartesian products and taking polar duals.
\end{dfn}
Hanner polytopes are presently the conjectured minimizers for \emph{the Mahler conjecture} (cf. \cite{Kim}\label{citas_20}).

\section{A relationship between {the}  matrix ${M_d(\ell,a)}$ and  the polar dual of  the {origin--symmetric} parallelepiped arising from ${M_d(\ell,a)}$}\label{sec:ring_of_mat_revisited}
By Remark \ref{rem:parallelepipeds} and Theorem \ref{thm:vol_isobiselado}, there is a  whole family of $d$--parallelepipeds having the
same volume as $\II_d(\ell,a)$: two such ones are $Par(D_d(\ell-a,a))$ and $Par({M_d(\ell,a)})$
by (\ref{eqn:diagonalization}) and (\ref{eqn:diagonal}). The former is a
$d$--box (cf. Notations \ref{nota:common}).   By Remark \ref{rem:parallelepipeds_2},  the latter is invariant under the
rotation $\sigma_{2\pi/d}:\R^d \rightarrow\R^d$. A translation $\tau:\R^d \rightarrow\R^d$ can be found  so that the origin $O$
becomes the barycenter of the image set $\tau\left(Par({M_d(\ell,a)})\right)$.

\begin{prop}\label{prop:vertices}
The vertices of  the polar dual of the origin--symmetric parallelepiped $\tau\left(Par({M_d(\ell,a)})\right)$ are the columns of the matrices $\pm 2\left({{M_d(\ell,a)}}\right)^{-1}$.
\end{prop}

\begin{proof}
A $d$--parallelepiped has $2d$ facets { and $2^d$ vertices}, so its polar dual has $2d$ vertices { and $2^d$ facets}.
A { standard} computation { in linear algebra}
shows  that the supporting hyperplanes of $\tau\left(Par({M_d(\ell,a)})\right)$ have
equations
{
\begin{equation}\label{eqn:facets_of_PP'_new}
(\ell+(d-2)a)x_i-a\sum_{j\neq i}x_j=\pm\frac{(\ell-a)(\ell+(d-1)a)}{2},\quad i\in[d].
\end{equation}
Indeed, the proof requires to write down two matrices of sizes $d\times 2^d$: the columns of the first (resp. second) matrix  are the  coordinates of the  vertices of $Par(M_d(\ell,a))$ (resp. $\tau(Par(M_d(\ell,a)))$). Then, for certain choices of columns of the second matrix, one must substitute   in the corresponding  equation  in (\ref{eqn:facets_of_PP'_new}). We only give a sample for $d=3$.
The coordinates of the vertices of $Par(M_3(\ell,a))$ are  the columns of the matrix
\begin{equation}
\begin{pmatrix}
  0 & \ell & a & a & \ell+a & \ell+a & 2a & \ell+2a \\
  0 & a & \ell & a & \ell+a & 2a & \ell+a & \ell+2a  \\
  0 & a & a & \ell & 2a & \ell+a & \ell+a & \ell+2a
\end{pmatrix}
\end{equation}
The translation $\tau$ of vector $\frac{-\ell-2a}{2}(e_1+e_2+e_3)$ yields that
the coordinates of the vertices of $\tau(Par(M_3(\ell,a)))$ are  the columns of the matrix
\begin{equation}\label{eqn:matrix_tau_par}
\frac{1}{2}\begin{pmatrix}
  -\ell-2a & \ell-2a & -\ell & -\ell & \ell & \ell & -\ell+2a & \ell+2a \\
  -\ell-2a & -\ell & \ell-2a & -\ell & \ell & -\ell+2a & \ell & \ell+2a \\
  -\ell-2a & -\ell & -\ell & \ell-2a & -\ell+2a & \ell & \ell & \ell+2a
\end{pmatrix}
\end{equation}
Columns 1,2,3 and 5 of (\ref{eqn:matrix_tau_par}) are the vertices of a facet of $\tau(Par(M_3(\ell,a)))$ and they  satisfy the equation
\begin{equation}
(\ell+a)x_3-a(x_1+x_2)=-\frac{(\ell-a)(\ell+2a)}{2}.
\end{equation}
Similar computations for the rest of  facets of $\tau(Par(M_3(\ell,a)))$.

An expression equivalent to (\ref{eqn:facets_of_PP'_new}) is
}
\begin{equation}\label{eqn:facets_of_PP'}
\frac{\ell+(d-2)a}{(\ell-a)(\ell+(d-1)a)}2x_i+\frac{-a}{(\ell-a)(\ell+(d-1)a)}2\sum_{j\neq i}x_j=\pm1,\quad i\in[d].
\end{equation}

The coefficients, written in columns,  of the left--hand--side of (\ref{eqn:facets_of_PP'}) yield  the following two (opposite) matrices
\begin{equation}\label{eqn:pm_matrix}
\pm 2\left(\frac{1}{\ell-a} I_d+\frac{-a}{(\ell-a)(\ell+(d-1)a)}J_d\right)=\pm 2({M_d(\ell,a)})^{-1},
\end{equation}
the last equality due to (\ref{eqn:inverse}) with $\alpha=\frac{1}{\ell-a}$, $\beta=\frac{-a}{(\ell-a)(\ell+(d-1)a)}$, $\beta'=a$.
These columns are the coordinates of the vertices of the polar dual of $\tau\left(Par({M_d(\ell,a)})\right)$.
\end{proof}

We wonder about generalizations of Proposition \ref{prop:vertices}.
\begin{question} Let ${\bf k}$ be a { subfield of $\R$} and  $\FF\subseteq \MM_d({\bf k})$ be a family of 
matrices closed under inverses. 
Each regular $M\in \FF$ gives rise (in the { sense of Notation \ref{nota:par}}) to a  $d$--parallelepiped $Par(M)\subseteq \R^d$ with $O$ at one vertex. Consider
the translated parallelepiped  $\tau (Par(M))$  centered at the origin $O$. Find conditions on $\FF$ so that the vertices of
$\left(\tau (Par(M))\right)^\circ$ are the columns of $\pm 2(M)^{-1}$.
\end{question}


\section{Volume of the  polar dual of an isocanted cube}\label{sec:vol_polar_isocanted}
Consider the following 
parameters (depending on $0<a<\ell$), homogeneous of degree $-1$
\begin{equation}\label{eqn:b_c}
b:=\frac{1}{\ell -a}>0,\qquad c:=\frac{2}{\ell }>0
\end{equation}
with relation
\begin{equation}\label{eqn:relation}
a=\frac{2b-c}{bc}.
\end{equation}
Let $(e_1,e_2,\ldots,e_d)$ be the \emph{standard vector basis} in  $\R^d$.
\begin{nota}[Vectors  and points $m_{i,j}(b,c)$]\label{nota:mij}
Let
\begin{equation}-m_{j,i}(b,c)=m_{i,j}(b,c)=\begin{cases}
ce_i,& \text{if \ }{ i\in[d],}\  j={0},\\
b(e_i-e_j),& \text{if \ } i,j\in[d], i\neq j.
\end{cases}\end{equation}
\end{nota}
Notice that $m_{i,i}(b,c)$  is not defined. We have   $(d+1)d$ different vectors/points. We write $m_{i,j}$ whenever $b,c$ are understood. These are called \emph{molecules} in \cite{Alexander_al}.\label{citas_21}

\begin{lem}[Polar dual of the isocanted cube]\label{lem:polar_dual_of_iso}
\begin{equation}\label{eqn:polar_4}
\II^\circ_d(\ell,a):=\conv(m_{i,j}(b,c):i,j\in{[d]\cup\{0\}}, i\neq j).
\end{equation}
\end{lem}
{
\begin{proof}
Use  (\ref{eqn:isocanted_by_ineqs}),  (\ref{eqn:polar_3}) and (\ref{eqn:b_c}).
\end{proof}
}
Clearly, $\II^\circ_d(\ell,a)$ is origin--symmetric. We  will write
\begin{equation}\label{eqn:J}
\JJ_d(b,c):=\II^\circ_d(\ell,a)
\end{equation}
where (\ref{eqn:b_c}) holds.  We write $\JJ_d$ and $\II_d^\circ$ whenever $b,c$ and $\ell,a$ are understood .

This polytope is   \emph{sparse} (i.e., most coordinates of each vertex
are zero)  and \label{obs:sparse_q}  and   if $\ell,a$ are rational numbers, then  it is a
\emph{rational} polytope. 


Polarity 
transforms vertices into facets. In \cite{Puente_Claveria_Iso}\label{citas_22} the following is proved. Vertices of $\II_d$ are in bijection with proper subsets of ${[d]\cup\{0\}}$. Given $W\subset{[d]\cup\{0\}}$, if $\underline{W}$ denotes a vertex
in $\II_d$, then the corresponding facet in $\II^\circ_d$ is $\overline{W}:=\conv\left(m_{i,j}:i\in W, j\in{[d]\cup\{0\}}\setminus W \right)$.
In particular, if $w=|W|$ then $\overline{W}$ has $w(d+1-w)$ vertices. Furthermore, $|W|=1$ or $d$ if and only if $\overline{W}$ is
a simplex. The polytope  $\II_d$ is {not simple; it is } almost--simple \label{dfn:almost_simple} and cubical, (cf. Notations \ref{nota:common})  where vertices with $|W|=1$ or $d$
have valency $d$ and the remaining vertices have valency $d+1$. Besides $\underline{W}$ and $\underline{W'}$ are joined by an edge in $\II_d$ if and only if  $W\subseteq W'$ and $|W'|=|W|+1$.

It easily follows from \cite{Puente_Claveria_Iso}\label{citas_34} that an edge joins $m_{i,j}$ and $m_{k,l}$ in
$\II^\circ_d$ if and only if $i=k$ or $j=l$, for $d\ge2$.

\bigskip


By   elementary methods we can prove:

$$\vol(\JJ_2)= \frac{2c}{2!}(2b+c)$$
$$\vol(\JJ_3)= \frac{2c}{3!}(6b^2+3bc+c^2)$$
$$\vol(\JJ_4)=\frac{2c}{4!}( 20b^3+10b^2c+4bc^2+c^3).$$

In the rest of this section we compute the exact volume of $\JJ_d(b,c)$, for arbitrary $d\ge2$.

\medskip

The following Euclidean distances are easy   to compute:
\begin{equation}\label{eqn:distances}
\dd(m_{i,j}, m_{r,t})=
\begin{cases} 2c,&\text{\ if\ } (i,j)=(t,r) \text{\ and\ } 0\in\{i,j,r,t\},\\
2b,&\text{\ if\ } (i,j)=(t,r) \text{\ and\ } 0\notin\{i,j,r,t\},\\
\sqrt{2}c,&\text{\ if\ } i=r=0 \text{\ or\ } j=t=0,\\ \sqrt{2}b,&\text{\ if\ } i=r \text{\ or\ } j=t \text{\ and\ } 0\notin\{i,j,r,t\},\\
{\sqrt{b^2+(b-c)^2}},&\text{\ if\ } i=r\neq 0 \text{\ and\ } j  \neq t=0.\\
\end{cases}
\end{equation} 

In  \cite{Puente_Claveria_Iso}\label{citas_23}  we computed the $f$--vector of $\II_d$ (which does not depend on $\ell,a$).\footnote{ Recently we noticed that the same vector appears in formula (28) in \cite{Shephard}\label{citas_24}.} Reversing this vector,  we obtain the  $f$--vector (i.e., numbers $f_0{^d}, f_1^d,\ldots, f_{d-1}^d$ of  vertices, edges, etc.,
facets) of $\JJ_d$, namely
\begin{equation}\label{eqn:f_vector_polar_of_iso}
f_k^d=\left(2^{k+2}-2 \right){{d+1}\choose{d-1-k} }=\left(2^{k+2}-2 \right){{d+1}\choose{k+2} }, \quad k=0,1,\ldots,d-1.
\end{equation}
Superscripts in (\ref{eqn:f_vector_polar_of_iso}) will be omitted, when understood.

{ $\JJ_d$ is not combinatorially equivalent to the type A root polytope $\conv(e_i-e_j: i,j\in [d+1])$ since the number of vertices in $\JJ_d$ doubles the number of vertices in the root polytope. Indeed, $f_0^d=(2^2-2){{d+1}\choose{2}}$, using (\ref{eqn:f_vector_polar_of_iso}).
}

\begin{rem}\label{rem:M}
$\II_d^\circ(\ell,a)$ is an example of the unit ball in the Lipschitz--free space arising from a finite pointed  metric space $M$ (cf. \cite{Alexander_al,Godard}\label{citas_25}). Namely, consider  $(M={[d]\cup\{0\}},\dd)$  with marked point $0$ and
\begin{equation}\label{eqn:distance_metric_space}
\dd(i,j)=\begin{cases}
           0, & \mbox{if } i=j \\
           \ell-a, & \mbox{if } i,j\in[d], i\neq j \\
           \ell/2, & \mbox{otherwise}.
         \end{cases}
\end{equation}
The space of Lipschitz  functions $\phi:M\rightarrow \R$ with $\phi(0)=0$ is denoted $Lip_0(M)$. A  canonical choice of
norm makes
$Lip_0(M)$ a Banach space, called \emph{Lipschitz dual of $M$}. The unit closed ball in $Lip_0(M)$ is denoted $B_{Lip_0(M)}$.
There is a canonical predual, denoted $\FF(M)$ and called  \emph{Lipschitz--free space over $M$}. The closed unit ball in $\FF(M)$ is
denoted $B_{\FF(M)}$. 
We have
\begin{equation}\label{eqn:balls}
B_{\FF(M)}=\II_d^\circ(\ell,a), \qquad \II_d(\ell,a)=B_{Lip_0(M)}.
\end{equation}
$M$ becomes a graph with weights $\dd(i,j)$ on edges. As a weighted graph,  $M$ above is   the complete graph  with
$d+1$ vertices having, at most,  two different weights:  weight $\ell/2$ for edges meeting 0 and weight $\ell-a$ for the rest of edges.
\end{rem}


\subsection{Roofs}
\emph{Roofs}   generalize  isosceles triangles and isosceles trapezia\footnote{Trapezia are sometimes called trapezoids.} as well
as { right}  frusta and  prisms  { and also} products of two regular simplices.   Roofs are prismatoids, defined
below.\footnote{Slightly different definitions of prismatoid are found in the literature.} {Prismatoids and roofs have two
bases. To be precise, we should write \emph{right roofs},
because of the property  \ref{item:distance_roof} in Remark \ref{rem:roof}:  the distance between the barycenters of the two bases is equal to the roof height.}

\begin{dfn}[Prismatoid] \label{dfn:prismatoid} A $d$--\emph{prismatoid} is a $d$--polytope $P\subseteq \R^d$ such that there exist
two different parallel hyperplanes $H_1, H_2\subseteq \R^d$ the union of which contains all the vertices in $P$. After relabelling
the hyperplanes, we can assume that $P\cap H_1$ is a facet of $P$, which we denote $MB(P)$ and call \emph{major { base}}
of $P$. We denote $P\cap H_2$ by $mb(P)$ and call it \emph{minor { base}} or \emph{crest} of $P$.
{The hyperplanes can be renamed accordingly: $MH(P)$ and $mh(P)$.}
\end{dfn}

{It follows from the mere definition that
the three  sets of vertices are non--empty and satisfy
\begin{equation}\label{eqn:vertices}
\ver (P)=\ver (MB(P))\cup\ver (mb(P)),
\end{equation}
so that
\begin{equation}
P=\conv(MB(P)\cup mb(P)).
\end{equation}
The affine spans  of the bases of $P$ satisfy
\begin{equation}\label{eqn:aff_spans}
\aff(MB(P))={MH(P)},\qquad \aff(mb(P))\subseteq {mh(P)}
\end{equation}
and this inclusion  may be strict. Furthermore, $\dim MB(P)=d-1$ but  $\dim mb(P)\le d-1$ and the inequality may be strict. If $\dim mb(P)=\dim MB(P)$, then the roles of the bases can be exchanged.} We get
\begin{equation}\label{eqn:dim_prismatoid}
\dim P=\dim MB(P)+1.
\end{equation}

\begin{rem}\label{rem:prismatoids}
For some prismatoids, the pair $H_1,H_2$ is not unique (e.g., for a $d$--parallelepiped there are $d$
possible pairs). In such  cases, the choices are finitely many and we fix one pair. One must prove that the results do not depend on
the chosen pair.
\end{rem}

{
\begin{dfn}[Height of a prismatoid]\label{dfn:height}
$h(P):=\dd({MH(P)},{mh(P)})=\dd(\aff(MB(P)),\aff(mb(P)))$ is called  \emph{height} of the prismatoid $P$.
\end{dfn}
}

\begin{ex}\label{ex:prismatoids}
Some  examples of 3--prismatoids  are: tetrahedra (where $mb(P)$ is a point), frusta of pyramids (with any number of edges
 in the { base}) as well as  prisms (where $\dim mb(P)=\dim MB(P)=2$),
and some  so--called (in Architecture)
\emph{hipped roofs}
  (see figure \ref{fig_roof});  (here $\dim mb(P)=1$).
\end{ex}

\begin{nota}[Barycenters]
Let $P$ be a $d$--prismatoid, $G$ and $g$ be the  barycenters of $MB(P)$ and $mb(P)$ respectively. In general, $\dd(G,g)\neq h(P)$.
\end{nota}

Our key definition comes next.

{
\begin{dfn}[Roof]\label{dfn:right_roof}
For $d\ge1$, a $d$--prismatoid $P\subseteq \R^{d}$  is a \emph{roof} if there exist integers $0<C,V\le d$ and positive
real numbers $\ell_1,\ell_2,\ell_3$ such that
\begin{enumerate}
\item $d=V+C-1$, \label{item:dimension_roof} 
\item there exists an affine--Euclidean system of coordinates $y_1,y_2,\ldots,y_d$  in $\R^d$ such that the line passing through $G$ and $g$ is the $y_d$ axis and  the sections of $P$ by orthogonal  hyperplanes  are product of simplices: more
precisely \label{item:section_roof}
\begin{equation}\label{eqn:secc_P}
 \se(s,P) =\Delta_{V-1}\left(\ell_1\left(1-\frac{s}{h}\right) +\ell_2\frac{s}{h}\right)\times \Delta_{C-1}\left(\ell_1\left(1-\frac{s}{h} \right)\right), \qquad s\in [0,h]
\end{equation}
where $h=h(P)$ is the height of $P$ (cf. Definition \ref{dfn:height}).   Besides, $\secc(0,P)=\Delta_{V-1}(\ell_1)\times\Delta_{C-1}(\ell_1)=MB(P)$ and $\secc(h,P)=\Delta_{V-1}(\ell_2)=mb(P)$, (cf. Notation \ref{nota:Delta})

\item there exists a surjective map $\phi: \ver (MB(P))\rightarrow \ver (mb(P))$ which extends, by convex combinations, to a surjective map $\phi: MB(P)\rightarrow mb(P)$ such that $\phi(G)=g$ and has congruent fibers,\label{item:lay_under_map}
\item for each $v\in\ver (MB(P))$, the line segment $[v, \phi(v)]$ is an edge of $P$ and the length $\dd(v,\phi(v))$ is $\ell_3$, independent of $v$,\label{item:third_length_roof}

\end{enumerate}
\end{dfn}

We call $\phi$ the \emph{lay--under map of the roof} $P$.

\begin{rem}\label{rem:roof} For a roof $P$, we get
\begin{enumerate}
\item $\phi^{-1}(x)\simeq \Delta_{C-1}(\ell_1)$, for all $x\in mb(P)$ (congruent fibers),\label{item:fibers}
\item $\dd(G,g)=h(P)$. \label{item:distance_roof}
\end{enumerate}
\end{rem}

\begin{lem}With notations above, we have
\begin{equation}\label{eqn:ell_3}
\ell_3^2=h^2+\left(r_{C-1}(\ell_1)\right)^2+\left(r_{V-1}(\ell_1)-r_{V-1}(\ell_2)\right)^2,
\end{equation}
where $r_d(\ell)$
is the \emph{circumradius} of $\Delta_d(\ell)$. \footnote{In \cite{Kachanovich}, it is shown that $r_d(\ell)= \sqrt{\frac{d}{2(d+1)}}\ell$. Clearly, $r_0(\ell)=0$ and $r_1(\ell)=\ell/2$.\label{citas_35}}
\end{lem}

\begin{proof}
We know that $BM(P)=\Delta_{V-1}(\ell_1)\times \Delta_{C-1}(\ell_1)$ and $bm(P)=\Delta_{V-1}(\ell_2)$ and the dimension of $P$ is $d=V+C-1$.
By a rigid motion we can assume, without loss of generality, that the major  (resp. minor) hyperplane has equation $y_{d}=0$ (resp. $y_{d}=h$). Let $g_1\in\R^{V-1}$ (resp. $g_2\in\R^{C-1}$) denote $y$--coordinates of the barycenter of $\Delta_{V-1}(\ell_1)$,
(resp. $\Delta_{C-1}(\ell_1)$). By Notation \ref{nota:Delta}, we have $g_1=g_2=0$.  Then $G=O\in\R^{d}$ (resp. $g=(O,h)\in\R^{d-1}\times \R$) is the barycenter of  $BM(P)$ (resp. $bm(P)$). By \emph{Pythagoras theorem} applied
twice,  we have
\begin{equation}\label{eqn:ell_3_2}
\ell_3^2=t^2+h^2, \text{\ with\ } t^2=\left(r_{C-1}(\ell_1)\right)^2+\left(r_{V-1}(\ell_1)-r_{V-1}(\ell_2)\right)^2.
\end{equation}
\end{proof}

\begin{ex}[First example of prismatoid not roof] The half  $3$--cube
$P=\conv\left(p_1,p_2,p_3,p_4,p_5,p_6\right)$ with
$p_1=(0,0,0)$, $p_2=(2,0,0)$, $p_3=(2,2,0)$, $p_4=(0,2,0)$, $p_5=(0,2,2)$ and  $p_6=(2,2,2)$  is a prismatoid,
with hyperplanes  $MH(P): x_3=0$ and $mh(P): x_3=2$,  bases $MB(P)=\conv\left(p_1,p_2,p_3,p_4\right)\simeq\Delta_1(2)\times\Delta_1(2)$
a square and
$mb(P)\simeq\conv\left(p_5,p_6\right)\simeq\Delta_1(2)$ a segment,  $C=V=2$ and $\ell_1=\ell_2=2$. We have
$G=(1,1,0)$ and $g=(1,2,2)$.

The  natural choice for \emph{lay--under map} is $\phi(p_1)=\phi(p_4)=p_5$ and $\phi(p_2)=\phi(p_3)=p_6$. Property \ref{item:third_length_roof} in
Definition \ref{dfn:right_roof} is not satisfied because $\dd(p_4,p_5)=2<2\sqrt2=\dd(p_1,p_5)$.  In addition,  Property \ref{item:distance_roof} in
Remark \ref{rem:roof} is not satisfied because  $\dd(G,g)=\sqrt5\neq h(P)=2$. Other choices of $\phi$ do not work either
(there are finitely many).

\end{ex}

\begin{ex}[Second example of prismatoid not roof] The snub half $3$--cube
$P=\conv\left(p_1,p_2,p_3,p_4,p_5,p_6'\right)$ with $p_1,p_2,p_3,p_4,p_5$ as above,
and  $p_6'=(2,0,2)$  is a prismatoid,
with hyperplanes and $MB(P)$ as above,   and
$mb(P)\simeq\conv\left(p_5,p_6'\right)\simeq\Delta_1(2\sqrt2)$ a segment,  $C=V=2$, $\ell_1=2$ and $\ell_2=2\sqrt2$. We have
$G=(1,1,0)$.

The  natural choice for \emph{lay--under map} is $\phi(p_1)=\phi(p_4)=p_5$ and $\phi(p_2)=\phi(p_3)=p_6'$ and yields
$\phi(G)=(1,1,2)=g$. However, properties  \ref{item:third_length_roof} and \ref{item:section_roof} in Definition \ref{dfn:right_roof}
are not satisfied, as can be checked easily.
Property \ref{item:distance_roof}
in Remark \ref{rem:roof} is  satisfied because  $\dd(G,g)= h(P)=2$. Other choices of $\phi$ do not work either.

\end{ex}

}

\begin{nota}
A $d$--roof $P$ is determined (up to rigid motions) if we know the numbers $C,V,\ell_1,\ell_2$ and its height $h=h(P)$. We will  write $P{\simeq}\Roof(C,V,\ell_1,\ell_2,h)$.
\end{nota}

\begin{figure}[ht]
\centering
\includegraphics[width=6cm]{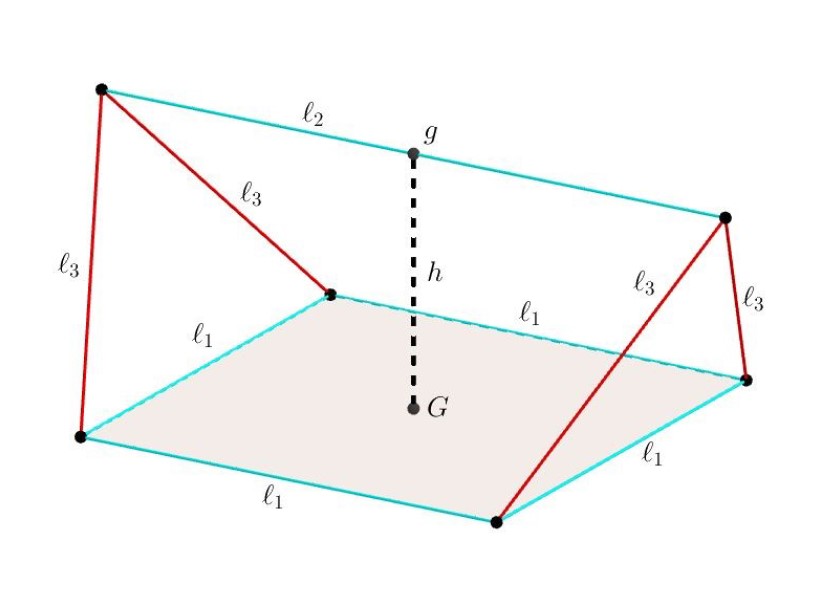}\ \includegraphics[width=6cm]{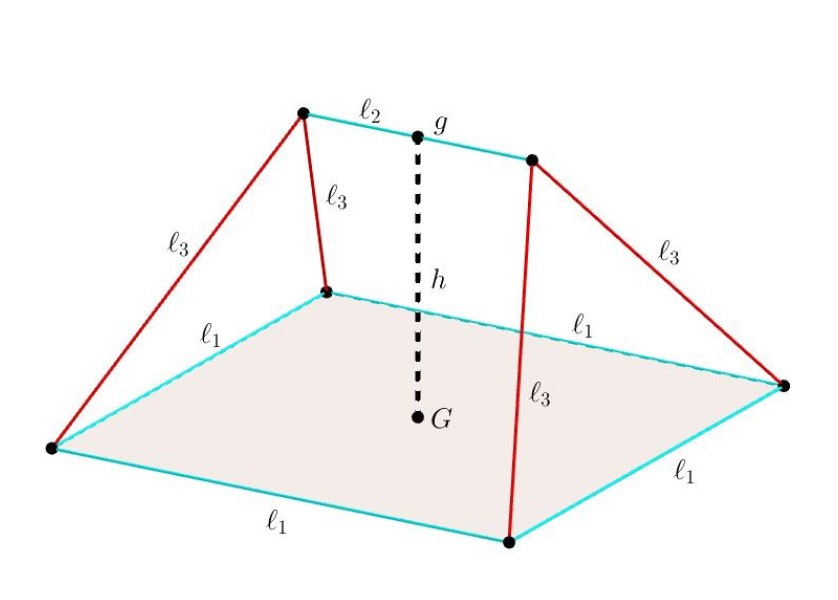}
\caption{Roofs   with $C=V=2$ and $d=3$. {Cases $\ell_1<\ell_2$ (left) and  $\ell_1>\ell_2$ (right).}}\label{fig_roof}
\end{figure}


\begin{rem}[Some particular cases of roofs]\label{rem:part_cases}
\begin{enumerate}
\item If   $V=1$, then $mb(P)$ is a point (and $\ell_2$ is irrelevant) and so $\Roof(C,1,\ell_1,\ell_2,h)$  is an isosceles $C$--simplex.
\item If   $C=1$, then the two bases $MB(P)$ and $mb(P)$ are regular simplices of  equal dimension.  Thus,  $\Roof(1,V,\ell_1,\ell_2,h)$  is a { right} frustum, if $\ell_1\neq\ell_2$, and a { right} prism, otherwise.

\end{enumerate}
\end{rem}


The volume of a $d$--prismatoid $P$ of height $h{=h(P)}$ is the integral
\begin{equation}
\vol_d P:=\int_{0}^h\vol_{d-1}\se(s,P)ds
\end{equation} where $\se(s,P)$ denotes the hyperplane section   of $P$ parallel to $MB(P)$ at level $s\in[0,h]$.
If $P$ is a roof then, { by definition,}  $\secc(s,P)$ is a
product of regular simplices.  Therefore, the exact volume of a roof can be obtained.

\begin{lem}[{Volume of roof section}]\label{lem:roof_section}
If $P{\simeq}\Roof(C,V,\ell_1,\ell_2,h)$ and   $d=V+C-1{=\dim P}$,
then for  each $s\in [0,h]$  it holds that
\begin{equation}\label{eqn:vol_section}
\vol_{d-1}\se(s,P)=\frac{1}{(C-1)!(V-1)!}\sqrt{\frac{CV}{2^{d-1}}}\sum_{n=0}^{V-1}{{V-1}\choose{n} }\ell_1^{C-1+n}\ell_2^{V-1-n}\left(1-\frac{s}{h} \right)^{C-1+n}\left(\frac{s}{h} \right)^{V-1-n}
\end{equation}
 is a  polynomial expression  in  $\frac{s}{h}$, homogeneous of degree  $d-1$.
\end{lem}
\begin{proof}
{ by definition, we have $\se(0,P)=MB(P)$, $\se(h,P)=mb(P)$. 
By (\ref{eqn:vol_product})  and  (\ref{eqn:known}),
we get
\begin{equation}
\vol_{d-1}\se(s,P)= \frac{1}{(V-1)!}\sqrt{\frac{V}{2^{V-1}}}\left(\ell_1 \left(1- \frac{s}{h}\right)+\ell_2 \frac{s}{h}\right)^{V-1} \frac{1}{(C-1)!}\sqrt{\frac{C}{2^{C-1}}}\ell_1^{C-1}\left( 1-\frac{s}{h}\right)^{C-1}
\end{equation}
}
 and we use the \emph{binomial expansion}
\begin{equation*}\label{eqn:binomial}
   \left(\ell_1 \left(1- \frac{s}{h}\right)+\ell_2 \frac{s}{h}\right)^{V-1}=\sum_{n=0}^{V-1}{{V-1}\choose{n} }\ell_1^n\left(1- \frac{s}{h}\right)^n\ell_2^{V-1-n}\left(\frac{s}{h} \right)^{V-1-n}
 \end{equation*}
 to conclude.
\end{proof}

\begin{lem}[Volume of roof]\label{lem:vol_roof}
If $P{\simeq}\Roof(C,V,\ell_1,\ell_2,{h})$ and    $d=V+C-1{=\dim P}$, then
\begin{equation}\label{eqn:vol_roof}
\vol_{d}P=\int_{0}^{h} { \vol_{d-1}}\se(s,P)ds=\frac{h}{d!}\sqrt{\frac{CV}{2^{d-1}}}\sum_{n=0}^{V-1}{{C-1+n}\choose{n} }\ell_1^{C-1+n}\ell_2^{V-1-n}.
\end{equation}
\begin{proof}
{Consider the number} ${k}(d,C,V):=\frac{1}{(C-1)!(V-1)!}\sqrt{\frac{CV}{2^{d-1}}}$.
From  (\ref{eqn:integrales}) { and (\ref{eqn:vol_section})} we get
\begin{align*}
\int_{0}^{h}{ \vol_{d-1}}\se(s,P)ds&={k}(d,C,V)\sum_{n=0}^{V-1}{{V-1}\choose{n} }\ell_1^{C-1+n}\ell_2^{V-1-n} h \beta(C+n,V-n)\\&={k}(d,C,V)\sum_{n=0}^{V-1}{{V-1}\choose{n} }\ell_1^{C-1+n}\ell_2^{V-1-n} h \frac{(C-1+n)!(V-1-n)!}{d!}
\end{align*}
 and {(\ref{eqn:vol_roof})} follows by  simplification of  factorials.
\end{proof}
\end{lem}

\begin{ex}
\begin{enumerate}
  \item Let $P$ be the frustum of an equilateral--triangular--based right pyramid, having bases of areas $a^2$ and $b^2$. Then $C=1$, $V=3$, $d=3$ and  $P{\simeq}\Roof(1,3,\ell_1,\ell_2,h)$ with $\ell_1=\frac{2a}{\sqrt[4]{3}}$, $\ell_2=\frac{2b}{\sqrt[4]{3}}$, and (\ref{eqn:vol_roof}) yields
      \begin{equation}\label{eqn:egypt_2}
        \vol_3 P=\frac{h}{3!}\sqrt{\frac{3}{4}}\left(\ell_1^2+\ell_1\ell_2+\ell_2^2\right)=
        \frac{h}{3!}\sqrt{\frac{3}{4}}\left(\frac{4a^2}{\sqrt3}+\frac{4ab}{\sqrt3}+\frac{4b^2}{\sqrt3}\right)=\frac{h}{3}\left(a^2+ab+b^2\right)
      \end{equation}
      recovering the Egyptian formula in (\ref{eqn:egypt}).
  \item Let $P$ be any roof shown in figure \ref{fig_roof}. Then $C=V=2$, $d=3$ and (\ref{eqn:vol_roof}) yields
      \begin{equation}\label{eqn:vol_roof_fig}
        \vol_3P=\frac{h}{3!}\left({{2}\choose{1}}\ell_1^2+{{1}\choose{0}}\ell_1\ell_2\right)=\frac{h}{6}\left(2\ell_1^2+\ell_1\ell_2\right).
      \end{equation}
\end{enumerate}
\end{ex}

\subsection{General case $d\ge2$}\label{subsec:general_case}
Recall that the values of $b,c>0$ are given in terms of $0<a<\ell$ in (\ref{eqn:b_c}).
By (\ref{eqn:f_vector_polar_of_iso}), there are $f_{d-1}=2^{d+1}-2$ facets in  $\JJ_d(b,c)$. They are contained in  the following hyperplanes. Two hyperplanes of equations
\begin{equation}\label{eqn:extraordinaryes}
H_{[d],\pm}: x_1+x_2+\cdots+x_d\pm c=0
\end{equation}
and, for each proper subset  $I\subsetneq [d]$ and each sign  ($+$ or $-$),  hyperplanes
\begin{equation}\label{eqn:no_extraordinaryes}
H_{I,\pm}: b\sum_{i\in I}x_i+(b-c)\sum_{j\in I^{\complement}}x_j\pm bc=0.
\end{equation}
 The  facet corresponding  to $I,+$  (i.e.,  the facet of $\JJ_d$ contained in  $H_{I,+}$) is
\begin{equation}\label{eqn:P_I+}
P_{I,+}:=\conv\left(m_{j,i}: i\in I, j\in I^\complement{\cup\{0\}} \right).
\end{equation}
Similarly
\begin{equation}\label{eqn:P_I-}
P_{I,-}:=\conv\left(m_{i,j}: i\in I, j\in I^\complement{\cup\{0\}} \right).
\end{equation}
The two  facets   contained in $H_{[d],\pm}$
are called \emph{{ extraordinary} facets } and are denoted $P_{[d],\pm}$:
\begin{equation}
P_{[d],+}:=\conv\left(m_{0,k}: k\in [d] \right)
\end{equation}
\begin{equation}
P_{[d],-}:=\conv\left(m_{k,0}: k\in [d] \right)
\end{equation}
The facets $P_{I,\pm}$ are called \emph{ordinary facets}.
Since
\begin{equation}\label{eqn:minus_sign}
P_{I,-}=-P_{I,+},\quad P_{[d],-}=-P_{[d],+},
\end{equation}    computations
are similar for  both signs. {Thus, Proposition \ref{prop:ord_facet_is_roof} below and its consequences are only stated for sign $+$.}

For each $1\le V\le d-1$ let us consider 
\begin{equation}\label{eqn:radical}
{\rad_d(b,c,V)}:=\sqrt{Vb^2+(d-V)(b-c)^2}.
\end{equation}

\begin{lem}[Distance to the origin]\label{lem:distancia_al_O}
For each dimension $d\ge2$, { parameters $b,c>0$,} proper subset $I\subsetneq [d]$  and sign  $\pm$,
the distance  from the hyperplane $H_{I,\pm}$ to the origin $O$ is
\begin{equation}\label{eqn:distance_to_O}
  \dd(H_{I,\pm},O)=\frac{bc}{\rad_d(b,c,|I|)}
\end{equation}
and it only depends on $d,b,c$ and $|I|$.
\end{lem}
\begin{proof}
Apply the Lagrange  distance formula from Section \ref{sec:background}.
\end{proof}

{

{
For the roofs $P$ we will meet below, there is a system of affine--Euclidean coordinates $y_1,y_2,\ldots,y_d$ in $\R^d$ such that the  lay--under map $\phi$ is the restriction to $MB(P)$ of the  composition of the translation $\tau:\R^d \rightarrow \R^d$ such that $\tau(G)=g$, the projection  $\rho: mh(P) \rightarrow mh(P)$ orthogonal onto $\aff(mb(P))$ and scaling by factor $\frac{\ell_2}{\ell_1}$, i.e.,
\begin{equation}\label{eqn:phi_composition}
\phi=\frac{\ell_2}{\ell_1}\rho\circ\tau|_{MB(P)}.
\end{equation}
}

\begin{prop}[Ordinary facets are roofs]\label{prop:ord_facet_is_roof}
For each proper subset  $I\subsetneq [d]$, the facet $P_{I,+}$ of $\JJ_d(b,c)$ is a roof and its height, denoted $h_{I,+}$   only depends on $d,b,c$ and $|I|$.
More precisely
\begin{equation}\label{eqn:bases}
MB(P_{I,+})=\conv\left( m_{j,i}: i\in I, j\in I^\complement \right)\qquad mb(P_{I,+})=\conv\left(m_{0,i}: i\in I \right)
\end{equation}
\begin{equation}\label{eqn:PI+}
{P_{I,+}=}\Roof\left(|I^\complement|, |I|, \sqrt2b,\sqrt2c,{h_{I,+}}\right),
\end{equation}
 and
\begin{equation}\label{eqn:height}
h_{I,+}=\frac{\rad_d(b,c,|I|)}{\sqrt{(d-|I|)|I|}}.
\end{equation}
\end{prop}
\begin{proof}
{ Write  $d'=d-1$. Recall $m_{0,i}= -ce_i$ and $m_{j,i}=b(e_j-e_i)$.
We have
\begin{equation}
P_{I,+}=\conv\left( m_{j,i}: i\in I, j\in I^\complement\cup\{0\} \right)
\end{equation}
with supporting hyperplane $H_{I,+}$ of equation
\begin{equation}
b\sum_{i\in I}x_i+(b-c)\sum_{j\in I^{\complement}}x_j+ bc=0.
\end{equation}
Inside $H_{I,+}$, consider the following parallel $(d-2)$--dimensional affine subspaces:
\begin{equation}
MH_{I,+}:
\sum_{i\in I}x_i+b=0,\qquad mh_{I,+}:\sum_{i\in I}x_i+c=0,
\end{equation}
with
\begin{equation*}\label{eqn:belong}
m_{j,i}\in MH_{I,+},\qquad m_{0,i}\in mh_{I,+}, \qquad i\in I, \ j\in I\complement
\end{equation*}
This proves  that $P_{I,+}$ is a $d'$--prismatoid.
}

Set $V=|I|$, $C=d-V=|I^\complement|$, $\ell_1=\sqrt2b$,  $\ell_2=\sqrt2c$, { and $\ell_3=\sqrt{b^2+(b-c)^2}$}. By the computed distances (\ref{eqn:distances}), we have
\begin{equation}
MB(P_{I,+}){\simeq}\Delta_{C-1}(\ell_1)\times \Delta_{V-1}(\ell_1),\qquad mb(P_{I,+}){\simeq}\Delta_{V-1}(\ell_2)
\end{equation}
\begin{equation}
P_{I,+}=\conv(MB(P_{I,+})\cup mb(P_{I,+}))
\end{equation}

with  $d'=V+C-1=d-1$, which shows Property \ref{item:dimension_roof} in Definition \ref{dfn:right_roof}.

The map $\phi$ acting on the vertices of $MB(P_{I,+})$ is
  \begin{equation}\label{eqn:phi}
    \phi(m_{j,i}):=m_{0,i}, \qquad i\in I, j\in I^\complement.
  \end{equation}
Then
\begin{equation}
  \phi(G)=\phi\left(\frac{1}{CV}\sum_{i\in I}\sum_{j\in I^\complement}m_{j,i}\right)=\frac{1}{CV}\sum_{i\in I}\sum_{j\in I^\complement}\phi\left(m_{j,i}\right)=\frac{1}{V}\sum_{i\in I}\phi\left(m_{0,i}\right)=g.
\end{equation}
Each $x\in mb(P_{I,+})$ satisfies $x=\sum_{i\in I}\lambda_i m_{0,i}$, for some $\lambda_i\in[0,1]$ with $\sum_{i\in I}\lambda_i=1$.
Thus $\phi^{-1}(x)=\conv\left(x_j: j\in I^\complement\right)$, with $x_j:=\sum_{i\in I}\lambda_i m_{j,i}$. Each fiber
$\phi^{-1}(x)$ is congruent to $\Delta_{C-1}(\sqrt2b)$, because $\dd(x_j,x_k)=\dd(m_{j,i},m_{k,i})=\sqrt2b$, by (\ref{eqn:distances}),
for all $j\neq k\in I^\complement$ and $i\in I$.
This shows Property \ref{item:lay_under_map} in Definition \ref{dfn:right_roof}.
Property \ref{item:third_length_roof} follows from (\ref{eqn:phi}) and the computed distances (\ref{eqn:distances}).


\end{proof}
}

\begin{ex}\label{ex:ell_3}
For $b,c>0$,  $d=4$, $I=\{1,2\}$, $I^\complement=\{3,4\}$, let us describe the facet $P_{I,+}$ as a {$d'$--roof with $d'=3$.
The supporting hyperplane is $H_{I,+}\simeq \R^3$ with equation
\begin{equation}
b(x_1+x_2)+(b-c)(x_3+x_4)+bc=0
\end{equation}
and
\begin{equation}
P_{I,+}=\conv(m_{3,1},m_{4,1},m_{3,2},m_{4,2},m_{0,1},m_{0,2}).
\end{equation}
Inside $H_{I,+}$, the equation $x_1+x_2+b=0$ defines $MH(P_{I,+})$ and the equation $x_1+x_2+c=0$ defines $mh(P_{I,+})$. With standard affine--Euclidean methods, we compute the height of  $P_{I,+}$ as a 3--prismatoid
\begin{equation}\label{eqn:hPI+}
h(P_{I,+})=\dd(MH(P_{I,+}),m_{0,1})=\dd(\pi(m_{0,1}),m_{0,1})=\frac{\sqrt{b^2+(b-c)^2}}{\sqrt2}=\frac{\rad_3(b,c,2)}{2}
\end{equation}
where $\pi(m_{0,1})=\frac{1}{2}(-b-c,c-b,b,b)$ is the orthogonal projection of $m_{0,1}$ onto $MH(P_{I,+})$.

The bases of $P_{I,+}$ are
\begin{equation}
mb(P_{I,+})=\conv(m_{0,1},m_{0,2}), \qquad MB(P_{I,+})=\conv(m_{3,1},m_{3,2},m_{4,1},m_{4,2})
\end{equation}
a segment and a square, respectively, and
\begin{equation}
P_{I,+}=\conv(MB(P_{I,+})\cup mb(P_{I,+})).
\end{equation}
}
We have $C=2=V$ and  $\ell_1=\sqrt2b$ and  $\ell_2=\sqrt2c$. We get $3=d'=V+C-1$, which  is condition \ref{item:dimension_roof} in
Definition \ref{dfn:right_roof}.
The map $\phi$ acting on the vertices of $MB(P_{I,+})$ is
  \begin{equation}
    \phi(m_{j,i}):=m_{0,i}, \qquad i\in \{1,2\}, j\in \{3,4\}.
  \end{equation}


{
For $s\in[0,h]$, the section  $\sec (s,P_{i,+})$ is equal to $\conv\left(q_{3,1}(s),q_{3,2}(s),q_{4,2}(s),q_{4,1}(s) \right)$,
where $q_{j,i}(s):=q_{j,i}(s):=m_{j,i}\left(1-\frac{s}{h}\right)+m_{0,i}\frac{s}{h}$, which yields Property \ref{item:section_roof}
of Definition \ref{dfn:right_roof}.

The length of the edge of $P_{I,+}$ joining the vertex  $v=m_{j,i}$ in $MB(P_{I,+})$ with the vertex $\phi(v)=m_{0,i}$ in
$mb(P_{I,+})$ 
is $\dd(v,\phi(v))=\dd(m_{j,i},m_{0,i})=\sqrt{b^2+(b-c)^2}=\ell_3$,
by (\ref{eqn:distances}), and it does not depend on $i\in I$ { or $j\in I^\complement$} (see  $\ell_3$ in figure
\ref{fig_roof}). This is Property \ref{item:third_length_roof} in Definition \ref{dfn:right_roof}.

Additionally, on can check that $\ell_3=\sqrt2 h$ agrees with relations (\ref{eqn:ell_3}) and (\ref{eqn:hPI+}.}
\end{ex}



 \begin{lem}[Pyramids on  extraordinary facets]\label{lem:sobre_fac_extra}
 \begin{equation}
 \vol_d \pyr(P_{[d],\pm},O)=\frac{c^d}{d!}
 \end{equation}
 \end{lem}
 \begin{proof} Without loss of generality we may assume that the sign  is $+$. By (\ref{eqn:distances}),  we get
 $P_{[d],+}=\Delta_{{d-1}}(\sqrt2c)$, hence, by (\ref{eqn:known}), we get
 $\vol_{d-1} P_{[d],+}=\frac{1}{(d-1)!}\sqrt d c^{d-1}$. By (\ref{eqn:lagrange})
 we get  $\dd(H_{[d],+},O)=\frac{c}{\sqrt d}$ and then  (\ref{eqn:vol_pyramid}) yields the result.
 \end{proof}
\begin{lem}[Pyramids on ordinary facets]\label{lem:sobre_fac_ord}
 Let  $I\subsetneq [d]$ be a proper subset. Then
 \begin{equation}
 \vol_d\pyr (P_{I,\pm},O)=\frac{1}{d!}\sum_{n=0}^{V-1}{{d-V-1+n}\choose{n} } b^{d-V+n}c^{V-n}
 \end{equation}
 and only depends  on  $d,b,c$ and  $V=|I|$.
 \end{lem}
 \begin{proof}
 Without loss of generality, we may assume that  the sign  is $+$. By  Lemma \ref{lem:vol_roof} and Proposition \ref{prop:ord_facet_is_roof},
 with $C=d-V$, $\ell_1=\sqrt2b$, $\ell_2=\sqrt2c$
we get
 \begin{equation}
 \vol_{d-1}(P_{I,+})=\frac{h_{I,+}}{(d-1)!}\sqrt{(d-V)V}\sum_{n=0}^{V-1}{{d-V-1+n}\choose{n} }b^{d-V-1+n}c^{V-1-n}
 \end{equation}
 By { (\ref{eqn:distance_to_O})} and (\ref{eqn:vol_pyramid}) we get  the result.
 \end{proof}

Now we are ready to compute $\vol_d \JJ_d(b,c)$.\label{com:ready} Notice that the cancelation of all radicals coming from
(\ref{eqn:distance_to_O}) and (\ref{eqn:height}), together with the cancelation of all $\sqrt2$ coming
from   (\ref{eqn:known}) and (\ref{eqn:distances}), provides the polynomial
expression (\ref{eqn:vol_polar_iso}) with rational coefficients. {
The fact that $\vol_d\JJ_d(b,c)$ is a polynomial in the variables $b,c$ homogeneous of degree $d$ can already be seen from the fact that each vertex coordinate of $\JJ_d(b,c)$ is linear on $b$ and $c$ and, by choosing an adequate triangulation, we can compute the volume as a sum of determinants of matrices whose columns are convex combinations of vertices of $\JJ_d(b,c)$.
}

\begin{thm}\label{thm:vol_polar_iso}
 For each dimension $d\ge2$, { and parameters $b,c>0$} we have
\begin{equation}\label{eqn:vol_polar_iso}
\vol_d \JJ_d(b,c)= \frac{2}{d!}\sum_{j=0}^{d-1}{{d+j-1}\choose{j}}b^jc^{d-j},
\end{equation}
a polynomial expression in $b,c$, homogeneous of degree $d$, without term in
$b^d$.
\end{thm}
\begin{proof}
We know that  the origin $O$ belongs to the  interior of $\JJ_d$.
The volume of $\JJ_d$ is the \emph{sum of volumes of pyramids} $\pyr(P,O)$ on all  facets $P$
 of $\JJ_d$.
By Lemmas \ref{lem:sobre_fac_extra} and \ref{lem:sobre_fac_ord}, this sum is
\begin{equation}\label{eqn:gran_sum}
\vol_d (\JJ_d)=2\left(\frac{c^d}{d!}+\sum_{V=1}^{d-1} {{d}\choose{d-V}}\frac{1}{d!}\sum_{n=0}^{V-1}{{d-V-1+n}\choose{n}} b^{d-V+n}c^{V-n}\right)
\end{equation}
where the factor 2 is due to central symmetry of $\JJ_d$ and the factor ${d}\choose{d-V} $ is due to a choice of $V$ elements in  $[d]$, i.e., a choice of a proper subset   $I\subsetneq [d]$ with  $V=|I|$.

In the right--hand side  of (\ref{eqn:vol_polar_iso}) we separate  the  $c^d$ term, getting
\begin{equation}\label{eqn:vol_polar_iso_2}
\frac{2}{d!}\left(c^d+\sum_{j=1}^{d-1}{{d+j-1}\choose{j} }b^jc^{d-j} \right)
\end{equation}
For   $1\le j\le d-1$,  write  $V-n=d-j$,  and compute  the following sum, using (\ref{eqn:igualdad_combi}) with $r=j-1$:
\begin{equation}\label{eqn:comb_equal_2}
\sum_{V-n=d-j; 1\le V\le d-1; 0\le n\le V-1} {{d}\choose{d-V}}{{d-V-1+n}\choose{n}}={{d+(d-V-1+n)}\choose{(d-V)+n}}={{d+j-1}\choose{j}}
\end{equation}
and when we substitute (\ref{eqn:comb_equal_2})  in  (\ref{eqn:gran_sum}) we conclude, in view of (\ref{eqn:vol_polar_iso_2}),
that  (\ref{eqn:vol_polar_iso}) holds true.
\end{proof}

 Notice that $\frac{2}{\ell}=c=0$ is a root of (\ref{eqn:vol_polar_iso}), meaning that  $\II_d^\circ(\ell,a)$ shrinks to the origin $O$, when $\ell$ tends to $\infty$.

Using (\ref{eqn:b_c}), we get $\frac{b}{c}=\frac{\ell}{2(\ell-a)}$ and the following.
\begin{cor}
\label{cor:vol_polar_iso}
 For each dimension $d\ge2$ and parameters $0<a<\ell$, we have
\begin{equation}\label{eqn:vol_polar_iso_3}
\vol_d\II_d^\circ(\ell,a)= \frac{2^{d+1}}{\ell^d d!}\sum_{j=0}^{d-1}{{d+j-1}\choose{j}}\left(\frac{\ell}{2(\ell-a)}\right)^j
\end{equation}
a rational expression in $\ell,a$, homogeneous of degree $-d$.\qed
\end{cor}

{
\begin{rem}\label{rem:Alexander}
At the bottom of  p.31 in \cite{Alexander_al}, the volume of $B_{\FF(K_{d+1})}$ is computed to be $\frac{1}{d!}{{2d}\choose{d}}$. Passing to our notation, $B_{\FF(K_{d+1})}=\JJ_d(1,1/2)$. The formula of Alexander et al. is a particular case of (\ref{eqn:vol_polar_iso_3}). Indeed, this easily follows from the combinatorial equality (\ref{eqn:mi_combinatorial_eq}).
\end{rem}
}


\begin{rem}
In an earlier paper, we computed the volume of the  polar dual in dimension 3, in a more general setting. If one sets
$\ell=2$ and
$a=-x=-y=-z$ in expression  (21) from \cite{Puente_Claveria},\label{citas_26} one gets $\frac{a^2-7a+16}{3(2-a)}$, which  is equal
to $\vol_3 \II^\circ_3(2,a)$, as expected.
\end{rem}


\section{Proof of the Mahler conjecture for isocanted cubes}\label{sec:Mahler}


The volume product for $\II_d(\ell,a)$ is obtained from (\ref{eqn:r}) and (\ref{eqn:vol_polar_iso_3}).
Multiplying out, we get a \emph{Laurent polynomial} in $\ell,\ell-a$, and a standard polynomial  in  $x:=\frac{\ell-a}{\ell}$. Explicitly,
\begin{align}\label{eqn:prod_vols}
\vol_d(\II_d)\vol_d(\II_d^\circ)(x)&=\\
&\frac{2}{(d-1)!}\sum_{j=0}^{d-1}2^{d-j}{{d+j-1}\choose{j} }x^{d-j-1}+\label{eqn:first}\\
&\frac{2(1-d)}{d!}\sum_{j=0}^{d-1}2^{d-j}{{d+j-1}\choose{j} }x^{d-j}.
\end{align}
\emph{The Mahler conjecture} for isocanted cubes can be expressed as the following inequality:
 \begin{equation}
 \frac{4^d}{d!}\le \vol_d(\II_d)\vol_d(\II_d^\circ)(x), \quad \forall\ 0<x<1.
 \end{equation}
 Equivalently,
 \begin{equation}\label{eqn:Mahler_2}
 0\le p_d(x), \quad \forall\ 0<x<1,
 \end{equation}
 where
 \begin{equation}
 p_d(x):=d!\left(\vol_d(\II_d)\vol_d(\II_d^\circ)(x)-\frac{4^d}{d!}\right){=\sum_{k=0}^d a_kx^k}
 \end{equation}
{ with
 \begin{equation}\label{eqn:term_indep}
 a_0=4d{{2d-2}\choose {d-1}}-4^d,\quad a_d=-2^{d+1}(d-1)
 \end{equation}
 \begin{equation}\label{eqn:coef_k_simo}
 a_k= { 2^{k+1}}\left(2d{{2d-k-2}\choose {d-1}}-(d-1){{2d-k-1}\choose {d-1}} \right), \quad 1\le k\le d-1.
 \end{equation}
 \begin{proof}
  The cases $d=2$  and 3 are  trivially checked: $p_2(x)=-8x^2+8x$   { and $p_3(x)=-32(x-1)(x+\frac{1}{2})^2$ are  positive in the open interval $(0,1)$. Assume $d\ge { 3}$.}
 First, notice that 1 is a root of $p_d(x)$, because $x=1$ is equivalent to $a=0$ and this means that the polytope described
 in (\ref{eqn:isocanted_by_ineqs}) is just a $d$--cube.\footnote{Alternatively, one can show that
 $\sum_{k=0}^da_k=0$ by using the equality $\sum_{k=0}^n{{2n-k}\choose{n}}2^k=2^{2n}$ from (1.17) in \cite{Gould_II}\label{citas_27}.}
 Secondly, it is clear that the leading coefficient $a_d$ is negative. The independent term  $a_0$  is positive,  for $d\ge  { 3}$.
 Indeed, using (\ref{eqn:bound_central}),
 for $d  +\frac{1}{2}\ge \pi$
 we get $d>\sqrt{d^2- { \frac{1}{4}}}\ge\sqrt{\pi(d {  { -}\frac{1}{2}})}$ and therefore
 \begin{equation}\label{eqn:bound}
 4d{{2d-2}\choose {d-1}} { \ge} 4d\frac{4^{d-1}}{\sqrt{\pi(d {  { -}\frac{1}{2}})}}>4^d.
 \end{equation}
 Next, simplifying factorials  and solving for $k$  { with $1\le k\le d-1$}, we easily show that $0\le a_k$ if and only if $ { 1\le }k\le\frac{3d-1}{d+1}$.

 Altogether, for each dimension  $d\ge3$ there exists an integer $k(d):=\lfloor\frac{3d-1}{d+1}\rfloor$ such that $0\le a_k$
 if   $0\le k\le k(d)$ and $a_k< 0$ if $1+k(d)\le k\le d$. By \emph{Descartes rule of signs} (cf. \cite{Kostrikin} p. 312\label{citas_28}),
 the polynomial $p_d(x)$ has exactly one
 positive root. We know that $1$ is a root and $a_0=p_d(0)>0$, whence we conclude by continuity.
 \end{proof}
 }

 {
 \begin{rem} \label{rem:M_2}
 The metric space $M$ in Remark \ref{rem:M} satisfies the \emph{four--points condition of Buneman} (cf. \cite{Godard}\label{citas_29}). Indeed, take $p,q,r,s\in M$ pairwise different.
\begin{itemize}
  \item if $0\in\{p,q,r,s\}$, then $\dd(p,q)+\dd(r,s)=\dd(p,r)+\dd(q,s)=\dd(p,s)+\dd(q,r)=\frac{3\ell}{2}-a$.
  \item otherwise $\dd(p,q)+\dd(r,s)=\dd(p,r)+\dd(q,s)=\dd(p,s)+\dd(q,r)=2\ell-2a$.
\end{itemize}
From Theorem 4.2 in \cite{Godard} it follows that $\FF(M)$ is
isometric to a subspace (i.e., a central section) of an $L_1$--space.\footnote{{An arbitrary 3--dimensional alcoved polytope
(more general than an isocanted cube) appears in the proof of Theorem 4.2 in \cite{Godard}.}}
Since $M$ has $d+1$
elements, then $\FF(M)$ is actually isometric to a subspace
of $\ell_1^d$.\footnote{We can provide an explicit description  of the norm in this subspace.} This fully agrees with  the
equivalence proved
in \cite{Bolker}: $Z$ is a zonoid if and only if $Z^\circ$ is a central section of the unit ball of $L_1([0,1])$.
As we know,  $\II_d(\ell,a)$ is a  zonotope and (\ref{eqn:balls}) holds.

The volume product for $B_{Lip_0(M)}$ with $M$ complete graph on $d+1$ vertices and just  one weight  is computed in
Claim 5 \cite{Alexander_al}. It applies to   $B_{Lip_0(M)}=\II_d(\ell,\ell/2)$, i.e., only  when  $a=\ell/2$.

Minimality of the volume product of the unit ball is studied in \cite{Alexander_al} for fixed $d=|M|-1$, under the additional hypothesis that $B_{\FF(M)}$ is simplicial. This does not apply to our polytopes since
$B_{\FF(M)}=\II_d^\circ(\ell,a)$ is not simplicial, because $\II_d(\ell,a)$ is almost--simple,  but not simple; (cf. paragraph after (\ref{eqn:J}) and Subsection \ref{subsec:general_case}.)
\end{rem}

 }

\section{Acknowledgements}
We deeply thank two anonymous referees for posing  queries which helped improve clarity
and for going beyond trends. We thank J. L\'{o}pez--Abad for his listening and questioning.


\end{document}